\documentclass{amsart}
\usepackage{amsfonts}
\usepackage{pstricks,pst-plot,pst-node,pst-tree,amssymb}
\usepackage{pstricks-add}
\usepackage{graphicx}
\usepackage{caption}
\usepackage{mathrsfs}

\numberwithin{equation}{section}
\input{tcilatex}

\begin{document}
\title[constant sigma-k curvature]{Umbilic hypersurfaces of constant sigma-k
curvature in the Heisenberg group}
\author{Jih-Hsin Cheng}
\address{Institute of Mathematics, Academia Sinica, Taipei and National
Center for Theoretical Sciences, Taiwan, R.O.C.}
\email{cheng@math.sinica.edu.tw}
\thanks{}
\author{Hung-Lin Chiu}
\address{Department of Mathematics, National Central University, Chung Li,
32054, Taiwan, R.O.C.}
\email{hlchiu@math.ncu.edu.tw}
\urladdr{}
\author{Jenn-Fang Hwang}
\address{Institute of Mathematics, Academia Sinica, Taipei, Taiwan, R.O.C.}
\email{majfh@math.sinica.edu.tw}
\thanks{}
\author{Paul Yang}
\address{Department of Mathematics, Princeton University, Princeton, NJ
08544, U.S.A.}
\email{yang@Math.Princeton.EDU}
\urladdr{}
\subjclass{Primary: 35L80; Secondary: 35J70, 32V20, 53A10, 49Q10.}
\keywords{sigma-k curvature, Heisenberg group}
\thanks{}

\begin{abstract}
We study immersed, connected, umbilic hypersurfaces in the Heisenberg group $%
H_{n}$ with $n$ $\geq $ $2.$ We show that such a hypersurface, if closed,
must be rotationally invariant up to a Heisenberg translation. Moreover, we
prove that, among others, Pansu spheres are the only such spheres with
positive constant sigma-k curvature up to Heisenberg translations.
\end{abstract}

\maketitle




\bigskip



\section{\textbf{Introduction and statement of the results}}

In \cite{CCHY} we studied umbilic hypersurfaces in the Heisenberg group $%
H_{n}.$ We proved that Pansu spheres are the only umbilic spheres with
positive constant $p$(or horizontal)-mean curvature in $H_{n}$ up to
Heisenberg translations. In this paper, we want to extend the result to the
case of constant horizontal sigma-i curvature.

First let us define horizontal sigma-i curvature of a hypersurface $\Sigma $
in $H_{n}.$ Throughout this paper, we always assume $\Sigma $ is immersed
(say, $C^{2}$ smooth and $C^{\infty }$ smooth in the region of regular
points, see the definition of regular point below) and $n\geq 2.$ Let $\xi $
($J,$ resp.) denote the standard contact ($CR$, resp.) structure on $H_{n},$
defined by the kernel of the contact form%
\begin{equation*}
\Theta =dt+\sum_{j=1}^{n}(x_{j}dy_{j}-y_{j}dx_{j})
\end{equation*}

\noindent where $x_{1},..,x_{n},y_{1},..,y_{n},t$ are coordinates of $H_{n}.$
A point $p$ $\in $ $\Sigma $ is called singular if $\xi $ $=$ $T\Sigma $ at $%
p.$ Otherwise $p$ is called regular or nonsingular. For a regular point, we
define $\xi ^{\prime }$ $\subset $ $\xi \cap T\Sigma $ by%
\begin{equation*}
\xi ^{\prime }=(\xi \cap T\Sigma )\cap J(\xi \cap T\Sigma ).
\end{equation*}

\noindent Let ($\xi ^{\prime }$)$^{\bot }$ denote the space of vectors in $%
\xi $, perpendicular to $\xi ^{\prime }$ with respect to the Levi metric $G$
:$=$ $\frac{1}{2}d\Theta (\cdot ,J\cdot )$ $=$ $%
\sum_{j=1}^{n}[(dx_{j})^{2}+(dy_{j})^{2}].$ It is not hard to see $\dim (\xi 
$ $\cap $ $T\Sigma )$ $\cap $ ($\xi ^{\prime }$)$^{\bot }$ $=$ $1.$ Take $%
e_{n}$ $\in $ ($\xi $ $\cap $ $T\Sigma )$ $\cap $ ($\xi ^{\prime }$)$^{\bot
} $ of unit length. Define the horizontal normal $e_{2n}$ :$=$ $Je_{n}$. Let 
$\nabla $ denote the pseudohermitian connection associated to $(J,\Theta ).$
Let $\alpha $ be the function on $\Sigma $, such that $\alpha e_{2n}+T$ $\in 
$ $T\Sigma $ where $T$ :$=$ $\frac{\partial }{\partial t}.$ Define $%
J^{\prime }$ on $\xi \cap T\Sigma $ by%
\begin{equation*}
J^{\prime }=J\text{ on }\xi ^{\prime }\text{ and }J^{\prime }e_{n}=0.
\end{equation*}

\noindent We now have a symmetric shape operator $\mathfrak{S}$ $:$ $\xi
\cap T\Sigma \rightarrow \xi \cap T\Sigma ,$ defined by%
\begin{equation}
\mathfrak{S(}v\mathfrak{)}=-\nabla _{v}e_{2n}+\alpha J^{\prime }v.
\label{1-1}
\end{equation}

\noindent (see Proposition 2.2 in \cite{CCHY}). A regular point $p$ is
called an umbilic point if $\mathfrak{S}(\xi ^{\prime })$ $\subset $ $\xi
^{\prime }$ at $p$ and all the eigenvalues of $\mathfrak{S}$ acting on $\xi
^{\prime }$ are the same (see also Definition 2.1 in Section 2). Let us
denote this common eigenvalue by $k.$ On the other hand, by Proposition 2.3
in \cite{CCHY}, we have%
\begin{equation*}
\mathfrak{S}(e_{n})=le_{n}
\end{equation*}

\noindent for some real number $l.$ The sigma-m curvature is the m-th
symmetric function of eigenvalues of the shape operator $\mathfrak{S}$. For
instance, the sigma-1 curvature is nothing but $p$-(or horizontal) mean
curvature $H$, which is the trace of $\mathfrak{S.}$ In terms of $k$ and $l,$
we have%
\begin{equation*}
H=l+(2n-2)k
\end{equation*}

\noindent at an umbilic point. Similarly the sigma-i curvature, denoted as $%
\sigma _{i,n},$ reads%
\begin{equation}
\sigma _{i,n}=\binom{2n-2}{i-1}lk^{i-1}+\binom{2n-2}{i}k^{i}  \label{1.1}
\end{equation}

\noindent for $1\leq i\leq 2n-2$ while%
\begin{equation}
\sigma _{2n-1,n}=lk^{2n-2}.  \label{1.2}
\end{equation}

In \cite{CCHY}, we study the case $i$ $=$ $1,$ i.e., $\sigma _{1,n}$ $:=H$
is a positive constant. Among others, we show that any umbilic sphere of
positive constant $H$ in $H_{n},$ $n$ $\geq $ $2,$ is a Pansu sphere up to a
Heisenberg translation. Let us recall what Pansu spheres are. For any $%
\lambda $\TEXTsymbol{>}$0$, the Pansu sphere $S_{\lambda }$ is the union of
the graphs of the functions $f$ and $-f$, where%
\begin{equation}
f(z)=\frac{1}{2\lambda ^{2}}\left( \lambda |z|\sqrt{1-\lambda ^{2}|z|^{2}}%
+\cos ^{-1}{\lambda |z|}\right) ,\ \ \ |z|\leq \frac{1}{\lambda }.
\label{1.6}
\end{equation}

\noindent It~is~known~that $S_{\lambda }$ has $p$-(or horizontal) mean
curvature $H$ $=$ $2n\lambda $ (see Section 2 in \cite{CCHY} for instance).

In this paper, we will study umbilic hypersurfaces of constant general
sigma-i curvature $\sigma _{i,n}.$ We have the following result.

\bigskip

\textbf{Theorem A}. \textit{Let }$\Sigma $\textit{\ be an immersed,
connected, orientable, closed, umbilic hypersurface of }$H_{n},$\textit{\ }$%
n $\textit{\ }$\geq $\textit{\ }$2$\textit{\ with nonvanishing Euler number}$%
.$\textit{\ For a given }$i,$\textit{\ }$1\leq i\leq 2n-1,$ \textit{suppose }%
$\sigma _{i,n}$\textit{\ of }$\Sigma $\textit{\ is a positive constant. Then 
}$\Sigma $\textit{\ must be a Pansu sphere up to a Heisenberg translation.}

\bigskip

\textbf{Theorem B}. \textit{Let }$\Sigma $\textit{\ be an immersed,
connected, umbilic hypersurface in }$H_{n}$\textit{, }$n$\textit{\ }$\geq $%
\textit{\ }$2$\textit{. Let }$S_{\Sigma }$ \textit{denote the set of all
singular points in} $\Sigma .$ \textit{Then }

\textit{(a)} \textit{either }$\alpha ^{2}+k^{2}$ $\equiv $ $0$ \textit{on }$%
\Sigma $ \textit{or }$\alpha ^{2}+k^{2}$\textit{\ }$>$\textit{\ }$0$\textit{%
\ at all points in }$\Sigma \backslash S_{\Sigma }$\textit{.}

\textit{(b) Suppose }$\alpha ^{2}+k^{2}$\textit{\ }$\equiv $\textit{\ }$0$%
\textit{\ on }$\Sigma .$\textit{\ Then }$\Sigma $\textit{\ is congruent with
part of a hypersurface }$C^{n-1}\times \gamma \times R$ \textit{where }$%
\gamma $\textit{\ is a curve in the Euclidean plane }$C$\textit{\ with
signed curvature }$l$\textit{\ (}$=H$\textit{)}$.$

\textit{(c) Suppose }$\alpha ^{2}+k^{2}$\textit{\ }$>$\textit{\ }$0$\textit{%
\ at all points in }$\Sigma \backslash S_{\Sigma }$\textit{. Then }$\Sigma $%
\textit{\ is congruent with part of a rotationally invariant hypersurface.
Moreover, the radius of leaves in the associated foliation is }$\frac{1}{%
\sqrt{\alpha ^{2}+k^{2}}}.$

\bigskip

For the situation of $\sigma _{i,n}$ $=$ $c$ $\leq $ $0,$ we have the
following result.

\bigskip

\textbf{Theorem C}. \textit{Let }$\Sigma $\textit{\ be an immersed,
connected, umbilic hypersurface of }$H_{n},$\textit{\ }$n$\textit{\ }$\geq $%
\textit{\ }$2.$

\textit{(a) Suppose }$\sigma _{i,n}$\textit{\ }$=$\textit{\ }$c$\textit{\ }$%
= $\textit{\ }$0.$\textit{\ Then there is no closed such hypersurface }$%
\Sigma $.

\textit{(b) Suppose }$\sigma _{i,n}$\textit{\ }$=$\textit{\ }$c$\textit{\ }$%
< $\textit{\ }$0.$\textit{\ Then for }$i$ odd, $1\leq i\leq 2n-1,$ \textit{%
any closed such hypersurface }$\Sigma $\textit{\ must be a Pansu sphere up
to a Heisenberg translation.}

\bigskip

We remark that for the case $\sigma _{i,n}$\textit{\ }$=$\textit{\ }$c$%
\textit{\ }$<$\textit{\ }$0,$ $i$ $\geq $ $4$ and even, there may correspond
closed umbilic hypersurfaces which are not congruent with Pansu spheres (see
Figure 4.2 and Proposition 4.2). In Section 2 we give a proof of Theorem B.
In Section 3 we discuss solutions to an ODE for $(\alpha ,k)$ (see (\ref%
{2-1-0}))$.$ The main result is described in Theorem 3.1 and Theorem 3.3. In
Section 4, among others, we apply these ODE results Theorem 3.1 and Theorem
3.3 to prove Theorem A and Theorem C. In Appendix A we derive an "almost"
closed form of solutions to (\ref{2-1-0}).

\bigskip

\textbf{Acknowledgments}. J.-H. C. (P. Y., resp.) is grateful to Princeton
University (Academia Sinica in Taiwan, resp.) for the kind hospitality.
J.-H. C., H.-L. C., and J.-F. H. would like to thank the Ministry of Science
and Technology of Taiwan for the support of the following research projects:
MOST104-2115-M-001-011-MY2, MOST104-2115-M-008-003-MY2, and
MOST104-2115-M-001-009-MY2, resp.. P. Y. would like to thank the NSF of the
United States for the grant DMS 1509505. We would also like to thank
Ms.Yu-Tuan Lin and Ms.Li-Yu Huang for the computer assistance to draw
figures.

\bigskip

\section{Umbilicity and rotational symmetry}

Let $e_{\beta },$ $e_{n+\beta }$ $:=$ $Je_{\beta },$ $1$ $\leq $ $\beta $ $%
\leq $ $n,$ be an orthonormal frame with respect to the Levi metric $G$ :$=$ 
$\frac{1}{2}d\Theta (\cdot ,J\cdot )$ $=$ $%
\sum_{j=1}^{n}[(dx_{j})^{2}+(dy_{j})^{2}]$ (see Section 1)$.$ Recall that $T$
$:=$ $\frac{\partial }{\partial t}$ is the (characteristic or Reeb) vector
field satisfying $\Theta (T)$ $=$ $1$ and $d\Theta (T,\cdot )$ $=$ $0.$ Let $%
(p;e_{\beta }(p),e_{n+\beta }(p),T(p))$ be a moving frame, depending on $%
p\in H_{n}$. There exist one-forms $\omega ^{a},\omega ^{2n+1},\omega
_{a}{}^{b}$ such that 
\begin{equation}
\begin{split}
dp& =e_{\beta }\otimes \omega ^{\beta }+e_{n+\beta }\otimes \omega ^{n+\beta
}+T\otimes \omega ^{2n+1} \\
de_{\gamma }& =e_{\beta }\otimes \omega _{\gamma }{}^{\beta }+e_{n+\beta
}\otimes \omega _{\gamma }{}^{n+\beta }+T\otimes \omega ^{n+\gamma } \\
de_{n+\gamma }& =e_{\beta }\otimes \omega _{n+\gamma }{}^{\beta }+e_{n+\beta
}\otimes \omega _{n+\gamma }{}^{n+\beta }-T\otimes \omega ^{\gamma } \\
dT& =0,
\end{split}
\label{moequ}
\end{equation}%
where $\Theta =\omega ^{2n+1}$ and 
\begin{equation}
\begin{split}
\omega _{a}{}^{b}& =-\omega _{b}{}^{a},\ \ \text{for}\ 1\leq a,b\leq 2n, \\
\omega _{n+\alpha }{}^{n+\beta }& =\omega _{\alpha }{}^{\beta },\ \omega
_{\alpha }{}^{n+\beta }=-\omega _{n+\alpha }{}^{\beta },\ \ \ \text{for}\
1\leq \alpha ,\beta \leq n.
\end{split}%
\end{equation}%
The equations (\ref{moequ}) are called \textbf{the equations of motion}. Let 
$\Sigma $ be a hypersurface of $H_{n}.$ Recall that $H_{n}$ has a standard
contact structure $\xi $ and $\mathfrak{S}$, the symmetric shape operator
acting on $\xi \cap T\Sigma $ is defined by (\ref{1-1}):%
\begin{equation*}
\mathfrak{S(}v\mathfrak{)}=-\nabla _{v}e_{2n}+\alpha J^{\prime }v
\end{equation*}

\noindent for $v$ $\in $ $\xi \cap T\Sigma .$ Recall that $\xi ^{\prime }$ $%
:=$ $(\xi \cap T\Sigma )\cap J(\xi \cap T\Sigma ).$ Recall that the
coefficients $h_{ab}$ of the second fundamental form are defined by%
\begin{equation*}
h_{ab}:=-G(\nabla _{e_{b}}e_{2n},e_{a})
\end{equation*}

\noindent for an orthonormal basis $e_{1},$ $..,$ $e_{n},$ $e_{n+1}$ $:=$ $%
Je_{1},$ $..,$ $e_{2n-1}$ $:=$ $Je_{n-1}$ with respect to the Levi metric $%
G. $

\bigskip

\textbf{Definition 2.1.} Let $\Sigma $ be a hypersurface of $H_{n}.$ Suppose 
$p$ $\in $ $\Sigma $ is a regular point. We say that $p$ is an umbilic point
if

\begin{equation}
\begin{split}
& (\text{i})\ \mathfrak{S}(\xi ^{\prime })\subset \xi ^{\prime }; \\
& (\text{ii})\ \mathfrak{S}=kI\text{ on }\xi ^{\prime }\text{ for some
constant }k.
\end{split}%
\end{equation}

\noindent where $I$ denotes the identity map.

\bigskip

Let us review some useful known results in \cite{CCHY}.

\bigskip

\textbf{Proposition 2.2 }(Proposition 4.1 in \cite{CCHY}). Suppose $\Sigma $
is an umbilic hypersurface of $H_{n}$. If $p\in \Sigma $ is a singular
point, then it is isolated.

\bigskip

\textbf{Proposition 2.3}. (Proposition 4.2 in \cite{CCHY}). Suppose that $%
\Sigma $ is an umbilic hypersurface of $H_{n}$. Then, on the regular part of 
$\Sigma $, we have 
\begin{equation}
\begin{split}
ek& =el=e\alpha =0,\ \ \ \text{for all}\ e\in \xi ^{^{\prime }}, \\
e_{n}k& =(l-2k)\alpha , \\
e_{n}\alpha & =k^{2}-\alpha ^{2}-kl.
\end{split}
\label{intcon}
\end{equation}

\textbf{Proposition 2.4} (Proposition 4.3 in \cite{CCHY}). Suppose that $%
\Sigma $ is an umbilic hypersurface. Then the subbundle $\xi ^{\prime }$
generates a $(2n-1)$-dimensional foliation on the regular part under the Lie
bracket. In addition, the characteristic direction $e_{n}$ is always
transversal to each leaf.

\bigskip

\textbf{Definition 2.5}. Let $\Sigma $ be a hypersurface in $H_{n}$. We say
that $\Sigma $ is rotationally invariant if it is invariant under the group
of rotations in $R^{2n+1}$ about $t$-axis, the last coordinate axis.

\bigskip

Recall (see Section 2 in \cite{CCHY}) that the coefficients $h_{ab}$ of the
second fundamental form are defined by%
\begin{equation*}
h_{ab}:=-G(\nabla _{e_{b}}e_{2n},e_{a})
\end{equation*}

\noindent for an orthonormal basis $e_{1},$ $..,$ $e_{n},$ $e_{n+1}$ $:=$ $%
Je_{1},$ $..,$ $e_{2n-1}$ $:=$ $Je_{n-1}$ with respect to the Levi metric $%
G. $

\bigskip

\textbf{Theorem 2.6 }(Proposition 3.1 in \cite{CCHY}). If $\Sigma $ is
rotationally invariant, then it is umbilic.

\bigskip

\proof
(An argument different from the one given in the proof of Proposition 3.1 in 
\cite{CCHY}) Since $\Sigma $ is rotationally invariant, it is easy to see
that $\Sigma $ has the induced foliation $\Sigma =\cup _{t}S_{\rho (t)}$,
where each leaf 
\begin{equation*}
S_{\rho (t)}=\{(z,t)\ |\ |z|=\rho (t)\}\ ,\ \text{for some}\ \rho (t)>0
\end{equation*}%
is a sphere. Actually this foliation is just the one generated by $\xi
^{\prime }=(\xi \cap T\Sigma )\cap J(\xi \cap T\Sigma )$ under the Lie
bracket. Taking a frame $\{e_{\beta },e_{n+\beta }=Je_{\beta }|\ \beta
=1,\cdots ,n-1\}$ of $\xi ^{\prime }$, we have, for $1\leq \beta \leq n-1$, 
\begin{equation*}
\lbrack e_{\beta },e_{n+\beta }]=-2(T+\alpha e_{2n}+\frac{h_{\beta \beta
}+h_{(n+\beta )(n+\beta )}}{2}e_{n}),\ \text{mod}\ \xi ^{\prime }.
\end{equation*}%
Since $[e_{\beta },e_{n+\beta }]$ is tangent to a leaf, there exist a $k$
such that $\frac{h_{\beta \beta }+h_{(n+\beta )(n+\beta )}}{2}=k$, for all $%
\beta $ ($k$ is independent of $\beta )$. We proceed to compute the
following Lie brackets, mod $\xi ^{\prime }$, 
\begin{equation*}
\begin{split}
\lbrack e_{\beta },T+\alpha e_{2n}+ke_{n}]& =(e_{\beta }k-e_{n+\beta }\alpha
-kh_{(n+\beta )n}-\alpha h_{n\beta })e_{n} \\
\lbrack e_{n+\beta },T+\alpha e_{2n}+ke_{n}]& =(e_{\beta }\alpha +e_{n+\beta
}k+kh_{\beta n}-\alpha h_{n(n+\beta )})e_{n}.
\end{split}%
\end{equation*}%
Since $\xi ^{\prime }$ generates a foliation and $e_{n}$ transverses to each
leaf, we have 
\begin{equation*}
\begin{split}
0& =e_{\beta }k-e_{n+\beta }\alpha -kh_{(n+\beta )n}-\alpha h_{n\beta } \\
0& =e_{\beta }\alpha +e_{n+\beta }k+kh_{\beta n}-\alpha h_{n(n+\beta )}.
\end{split}%
\end{equation*}%
Since $\Sigma $ is rotationally invariant, both $k$ and $\alpha $ are
constants on each leaf, we get 
\begin{equation*}
\begin{split}
0& =kh_{(n+\beta )n}+\alpha h_{n\beta } \\
0& =kh_{\beta n}-\alpha h_{n(n+\beta )}.
\end{split}%
\end{equation*}%
This implies that $h_{n\beta }=h_{n(n+\beta )}=0$ , for $1\leq \beta \leq
n-1 $, which means that $e_{n}$ is an eigenvalue of the shape operator $%
-\nabla e_{2n}+\alpha J^{\prime }$. (or $(-\nabla e_{2n}+\alpha J^{\prime
})(\xi ^{\prime })\subset \xi ^{\prime }$). Now we can assume we take a
frame of $\xi ^{\prime }$ such that $e_{\beta },e_{n+\beta }$ are
eigenvectors. To complete the proof, we need to show that $h_{\beta \beta
}=h_{(n+\beta )(n+\beta )}=k$, for $1\leq \beta \leq n-1$. Let $\nu =\alpha
e_{2n}+ke_{n}\in (\xi ^{\prime \perp }$, we have shown $[e_{\beta
},e_{n+\beta }]=-2(T+\nu )$, mod $\xi ^{\prime }$, which implies that $T+\nu
\in TS_{\rho (t)}$. Thus $\nu $ is uniquely determined and, in the case $%
S_{\rho (t)}$, we have $-\frac{J(\nu )}{|\nu |^{2}}=x_{\gamma }\mathring{e}%
_{\gamma }+y_{\gamma }\mathring{e}_{n+\gamma }$, for all $(x_{\gamma
},y_{\gamma },t)\in S_{\rho (t)}$. We compute $\alpha \nu +kJ\nu =(\alpha
^{2}+k^{2})e_{2n}$, hence 
\begin{equation*}
\begin{split}
e_{2n}& =\frac{\alpha \nu +kJ\nu }{\alpha ^{2}+k^{2}}=-\left( \alpha \left( 
\frac{-\nu }{|\nu |^{2}}\right) +k\left( \frac{-J\nu }{|\nu |^{2}}\right)
\right) \\
& =-(\alpha y_{\gamma }+kx_{\gamma })\mathring{e}_{\gamma }-(-\alpha
x_{\gamma }+ky_{\gamma })\mathring{e}_{n+\gamma }.
\end{split}%
\end{equation*}%
Therefore, for each unit vector $e=a^{\beta }\mathring{e}_{\beta
}+a^{n+\beta }\mathring{e}_{n+\beta }\in \xi ^{\prime }$, 
\begin{equation}
-\nabla _{e}e_{2n}=[\alpha a^{n+\gamma }+ka^{\gamma }]\mathring{e}_{\gamma
}+[-\alpha a^{\gamma }+ka^{n+\gamma }]\mathring{e}_{n+\gamma },
\end{equation}%
hence 
\begin{equation}
\begin{split}
\big<-\nabla _{e}e_{2n},e\big>& =\alpha a^{\gamma }a^{n+\gamma }+k(a^{\gamma
})^{2}-\alpha a^{\gamma }a^{n+\gamma }+k(a^{n+\gamma })^{2} \\
& =k\left[ (a^{\gamma })^{2}+(a^{n+\gamma })^{2}\right] =k.
\end{split}%
\end{equation}%
This completes the proof.

\endproof%

\bigskip

Now we are ready to prove Theorem B stated in Section 1.

\bigskip

\proof
\textbf{(of Theorem B)} To prove (b), note that if $\alpha $ $\equiv $ $0$
and $k$ $\equiv $ $0$ on an umbilic hypersurface $\Sigma ,$ then $H$ $=$ $%
l+(2n-2)k$ $=$ $l$ is a constant on each leave of the foliation for $\Sigma $
in view of Proposition 4.2 in \cite{CCHY} (noting that $el$ $=$ $0$ for all $%
e$ $\in $ $\xi ^{\prime }$ and $\hat{e}_{2n}l$ $=$ $Tl$ $=$ $0)$. Moreover,
if $\alpha $ $\equiv $ $0,$ then $\Sigma $ is locally congruent with the
hypersurface $\Sigma ^{\ast }$ $\times $ $I$ where $\Sigma ^{\ast }$ is a
hypersurface of $C^{n}$ $=$ $R^{2n}$ (cf. Theorem 1.5 in \cite{CCHY}) and $I$
is an open interval of $R^{1}$.

Since $\alpha =k=0$, from \textbf{the equations of motion }(see (\ref{moequ}%
)), we have 
\begin{equation}
\begin{split}
de_{n}& =le_{2n}\otimes \omega ^{n}+T\otimes \omega ^{2n}, \\
de_{2n}& =-le_{n}\otimes \omega ^{n}-T\otimes \omega ^{n}.
\end{split}
\label{me1}
\end{equation}%
Here we have used $\omega _{2n}$ $^{a}(T)$ $=$ $\omega _{2n}$ $^{a}(\hat{e}%
_{2n})$ $=$ $0$ since $\alpha $ $=$ $0,$ and hence $\omega _{2n}$ $^{a}=$ $%
-l\delta _{n}^{a}\omega ^{n}$ (see Proposition 5.5 in \cite{ChiuLai}). Let $%
Z_{n}=\frac{1}{2}(e_{n}-ie_{2n}),\ \theta ^{n}=\omega ^{n}+i\omega ^{2n}$,
the complex version of (\ref{me1}) is 
\begin{equation}
dZ_{n}=ilZ_{n}\otimes \omega ^{n}+\frac{1}{2}iT\otimes \theta ^{\bar{n}}.
\label{me2}
\end{equation}%
Therefore, we have 
\begin{equation}
\begin{split}
d_{e_{n}}Z_{n}& =ilZ_{n},\ \ \text{mod}\ T; \\
d_{e}Z_{n}& =0,\ \ \text{mod}\ T,\ \text{for all}\ e\in \xi ^{\prime }.
\end{split}
\label{me3}
\end{equation}%
On the other hand, $\alpha =k=0$ implies that, for $1\leq \beta \leq n-1$, 
\begin{equation*}
\lbrack e_{\beta },e_{n+\beta }]=-2(T+\alpha e_{2n}+ke_{n})=-2T,\ \ \text{mod%
}\ \xi ^{\prime }.
\end{equation*}%
This means that $T$ is tangent to each leaf $L$ of the $(2n-1)$-dimensional
foliation of $\Sigma $ (see Proposition 2.4). Thus $Z_{n}\perp L$ with
respect to the adapted metric. Let $\pi :H_{n}\rightarrow C^{n}$ be the
projection along $T$, and $Z_{n}^{\ast }=\pi _{\ast }Z_{n},\ L^{\ast }=\pi
(L)$, we have $Z_{n}^{\ast }\perp L^{\ast }$ with respect to the standard
Euclidean metric on $C^{n}$. Here $L^{\ast }$ constitute a $(2n-2)$%
-dimensional foliation on $\Sigma ^{\ast }$. Also, from the second equation
of (\ref{me3}), we see that $Z_{n}^{\ast }$ is a constant along each leaf $%
L^{\ast }$.

Let $\Gamma(s)$ be a characteristic curve such that $\Gamma^{\prime
}(s)=e_{n}$ and let $\gamma(s)=\pi(\Gamma(s))$. We have $\gamma^{\prime
}(s)=\pi_{*}e_{n}=e_{n}^{*}$. We compute the derivative of $Z_{n}^{*}$ along 
$\gamma(s)$. From (\ref{me3}), we have 
\begin{equation}  \label{me4}
\begin{split}
Z_{n}^{*}{}^{\prime }(s)&=Z_{n}{}^{\prime }(s),\ \ \text{mod}\ T \\
&=d_{e_{n}}Z_{n},\ \ \text{mod}\ T \\
&=ilZ_{n}^{*}(s),
\end{split}%
\end{equation}
which implies that 
\begin{equation}  \label{me5}
Z_{n}^{*}(s)=e^{i[l(s)-l(s_{0})]}Z_{n}^{*}(s_{0}),\ \ \text{for a fixed}\
s_{0}.
\end{equation}
Therefore $\gamma$ completely lies on the complex line spanned by the
complex vector $Z_{n}^{*}(s_{0})$. After an unitary transformation, we can
assume, without loss of generality, that this complex line is spanned by the
last coordinates $z_{n}$. Since $Z_{n}^{*}(s_{0})\perp L^{*}(s_{0})$, we see
that $L^{*}(s_{0})$ is an open subset of $C^{n-1}$ with coordinates $%
(z_{1},\cdots,z_{n-1})$. Thus we get $\Sigma^{*}=L^{*}(s_{0})\times\gamma$
with $L^{*}(s_{0})\subset C^{n-1}$ an open subset, $\gamma\subset C^{1}$,
hence $\Sigma=L^{*}(s_{0})\times\gamma\times I$.

Finally, we compute 
\begin{equation}  \label{me6}
\begin{split}
l&=\big<\nabla_{e_n}e_{n},e_{2n}\big>=\big<d_{e_n}e_{n},e_{2n}\big> \\
&=\big<\frac{d}{ds}\Gamma^{\prime }(s),e_{2n}\big>=\frac{d}{ds}%
\gamma^{\prime }(s)\cdot \pi_{*}e_{2n},
\end{split}%
\end{equation}
where $\cdot$ is the standard inner product on the Euclidaen plane $C^{1}$.
Since $\pi_{*}e_{2n}=\hat{J}\gamma^{\prime }$, where $\hat{J}$ is the
complex structure on $C^{1}$, we see that it is normal to $\gamma$, thus,
from (\ref{me6}), $l$ is just \textbf{the signed curvature} of $\gamma$ with
respect to the orthonormal frame field $\{\pi_{*}e_{n},\pi_{*}e_{2n}\}$. We
have proved (b).

Next we are going to prove (c). On $\Sigma \backslash S_{\Sigma }$, we
define 
\begin{equation}
N=\frac{\alpha e_{n}-ke_{2n}}{\sqrt{\alpha ^{2}+k^{2}}}
\end{equation}%
(note that $\alpha ^{2}+k^{2}$ $>$ $0$ on $\Sigma \backslash S_{\Sigma }$ by
assumption). We will show that $N$ is the unit horizontal normal of each
leaf $L$ in the following sense 
\begin{equation}
d\left( X-\frac{N}{\sqrt{\alpha ^{2}+k^{2}}}(X)\right) \Big|_{L}=0
\label{lecircle}
\end{equation}%
where $X$ denotes the position vector. Since $L$ is spanned by $e_{\gamma
},e_{n+\gamma }$ and their Lie brackets, we see that (\ref{lecircle}) is
equivalent to 
\begin{equation}
\begin{split}
0& =d_{e_{\beta }}\left( X+\frac{-\alpha e_{n}+ke_{2n}}{\alpha ^{2}+k^{2}}%
(X)\right) \\
0& =d_{e_{n+\beta }}\left( X+\frac{-\alpha e_{n}+ke_{2n}}{\alpha ^{2}+k^{2}}%
(X)\right) .
\end{split}
\label{keycom}
\end{equation}%
Here $d_{v}X$ means differentiation of $X$ in the direction $v.$ Note that $%
d_{v}X$ $=$ $v.$ Now we are going to show (\ref{keycom}): By the motion
equation, for each $\gamma ,\ 1\leq \gamma \leq n-1$, we have 
\begin{equation}
\begin{split}
d_{e_{\gamma }}e_{n}& =\big<e_{\beta }\otimes \omega _{2n}{}^{n+\beta
}-e_{n+\beta }\otimes \omega _{2n}{}^{\beta }+T\otimes \omega
^{2n},e_{\gamma }\big> \\
& =-\sum_{\beta =1}^{n-1}h_{(n+\beta )\gamma }e_{\beta }+\sum_{\beta
=1}^{n}h_{\beta \gamma }e_{n+\beta } \\
& =-h_{(n+\gamma )\gamma }e_{\gamma }+h_{\gamma \gamma }e_{n+\gamma } \\
& =\alpha e_{\gamma }+ke_{n+\gamma }.
\end{split}
\label{keycom1}
\end{equation}%
Similarly we have 
\begin{equation}
\begin{split}
d_{e_{n+\gamma }}e_{n}& =-ke_{\gamma }+\alpha e_{n+\gamma } \\
d_{e_{\gamma }}e_{2n}& =-ke_{\gamma }+\alpha e_{n+\gamma } \\
d_{e_{n+\gamma }}e_{2n}& =\alpha e_{\gamma }-ke_{n+\gamma }.
\end{split}
\label{keycom2}
\end{equation}%
Recall that both $\alpha $ and $k$ are constants on each leaf, hence, using
formula (\ref{keycom1}), (\ref{keycom2}), we have (\ref{keycom}), hence (\ref%
{lecircle}). The formula (\ref{lecircle}) means that for each leaf $L$,
there exists a constant $X_{L}\in H_{n}$ such that%
\begin{equation}
X-X_{L}=\frac{N}{\sqrt{\alpha ^{2}+k^{2}}}(X),\ \text{for all}\ X\in L.
\label{keycom2-1}
\end{equation}%
We write $\frac{N}{\sqrt{\alpha ^{2}+k^{2}}}(X)=\sum_{j=1}^{2n}a_{j}(X)%
\mathring{e}_{j}(X).$ Here $\mathring{e}_{j}(X)$ $:=$ $\frac{\partial }{%
\partial x_{j}}+x_{j}^{\prime }T$ in which $x_{n+j}$ $=$ $y_{j},$ $%
x_{j}^{\prime }$ $=$ $y_{j},$ and $x_{n+j}^{\prime }$ $=$ $-x_{j}$ for $1$ $%
\leq $ $j$ $\leq $ $n$. Then we have%
\begin{eqnarray}
&&\sum_{j=1}^{2n}a_{j}(X)(x_{j}^{\prime }-(X_{L})_{j}^{\prime })
\label{keycom2-2} \\
&=&G(X-X_{L},J(X-X_{L}))=0  \notag
\end{eqnarray}

\noindent in view of (\ref{keycom2-1}), and hence%
\begin{eqnarray}
&&\sum_{j=1}^{2n}a_{j}(X)\mathring{e}_{j}(X)  \label{keycom2-3} \\
&=&\sum_{j=1}^{2n}a_{j}(X)(\frac{\partial }{\partial x_{j}}+x_{j}^{\prime }T)
\notag \\
&=&\sum_{j=1}^{2n}a_{j}(X)(\frac{\partial }{\partial x_{j}}%
+(X_{L})_{j}^{\prime }T)+\sum_{j=1}^{2n}a_{j}(X)(x_{j}^{\prime
}-(X_{L})_{j}^{\prime })T  \notag \\
&=&\sum_{j=1}^{2n}a_{j}(X)(\frac{\partial }{\partial x_{j}}%
+(X_{L})_{j}^{\prime }T)+0\text{ \ (by (\ref{keycom2-2}))}  \notag \\
&=&\sum_{j=1}^{2n}a_{j}(X)\mathring{e}_{j}(X_{L}),  \notag
\end{eqnarray}

\noindent meaning $L$ is the sphere with radius $\frac{1}{\sqrt{\alpha
^{2}+k^{2}}}$ and center $X_{L}$. Finally, we will show that 
\begin{equation}
-\frac{k}{\alpha ^{2}+k^{2}}T=d_{e_{n}}\left( X+\frac{-\alpha e_{n}+ke_{2n}}{%
\alpha ^{2}+k^{2}}(X)\right) =d_{e_{n}}(X_{L}).  \label{keycom3}
\end{equation}%
This equation says that the center $X_{L}$ lies on a line parallel to the $t$%
-axis. Therefore, $\Sigma \backslash S_{\Sigma }$ (and hence $\Sigma ,$
since $\Sigma \backslash S_{\Sigma }$ is open and dense in $\Sigma $ by
Theorem D in \cite{chy}$)$ is congruent with a rotationally invariant
hypersurface. For (\ref{keycom3}), we compute, as we just did for (\ref%
{keycom1}), 
\begin{equation}
\begin{split}
d_{e_{n}}e_{n}& =\big<e_{\beta }\otimes \omega _{n}{}^{\beta }-e_{n+\beta
}\otimes \omega _{2n}{}^{\beta }+T\otimes \omega ^{2n},e_{n}\big> \\
& =-\sum_{\beta =1}^{n-1}h_{(n+\beta )n}e_{\beta }+\sum_{\beta
=1}^{n}h_{\beta n}e_{n+\beta } \\
& =le_{2n},
\end{split}
\label{keycom4}
\end{equation}%
Similarly 
\begin{equation}
d_{e_{n}}e_{2n}=-le_{n}-T.  \label{keycom5}
\end{equation}%
And, using (\ref{intcon}), we have 
\begin{equation}
\begin{split}
e_{n}\left( \frac{\alpha }{\alpha ^{2}+k^{2}}\right) & =d_{e_{n}}\left( 
\frac{\alpha }{\alpha ^{2}+k^{2}}\right) =\frac{k^{2}+\alpha ^{2}-kl}{\alpha
^{2}+k^{2}}, \\
e_{n}\left( \frac{k}{\alpha ^{2}+k^{2}}\right) & =d_{e_{n}}\left( \frac{k}{%
\alpha ^{2}+k^{2}}\right) =\frac{l\alpha }{\alpha ^{2}+k^{2}}.
\end{split}
\label{keycom6}
\end{equation}%
Therefore, by the motion equation and formula (\ref{keycom4}),(\ref{keycom5}%
) and (\ref{keycom6}), we obtain (\ref{keycom3}). This completes the proof
of (c).

To prove (a), suppose the converse holds. Then there exists a sequence of $%
p_{j}$.$\in $.$\Sigma \backslash S_{\Sigma }$ with $\alpha _{j}^{2}$ $+$ $%
k_{j}^{2}$ $>$ $0,$ converging to $p_{\infty }$ $\in $ $\Sigma $ with $%
\alpha _{\infty }^{2}$ $+$ $k_{\infty }^{2}$ $=$ $0.$ Here we denote $\alpha
(p_{j})$ and $k(p_{j})$ by $\alpha _{j}$ and $k_{j},$ resp. for $j$ $=$ $1,$ 
$2,$ $..,\infty .$ Consider the vector

\begin{equation*}
v_{j}:=p_{j}-X_{L_{j}}
\end{equation*}%
\noindent where $L_{j}$ denotes the leaf through $p_{j}$ with the center $%
X_{L_{j}}$. From (\ref{keycom2-1}) and (\ref{keycom2-3}), we learn that $%
v_{j}$ is sitting in the contact plane passing through $X_{L_{j}}$ for all $%
j $ and 
\begin{equation}
||v_{j}||_{G}:=\sqrt{G(v_{j},v_{j})}=\frac{1}{\sqrt{\alpha _{j}^{2}+k_{j}^{2}%
}}  \label{vec-1}
\end{equation}

\noindent (recall that $G$ denotes the Levi metric). As $j$ $\rightarrow $ $%
\infty ,$ $||v_{j}||_{G}$ goes to infinity since $\sqrt{\alpha
_{j}^{2}+k_{j}^{2}}$ tends to zero in (\ref{vec-1}). So by (\ref{keycom3}), $%
X_{L_{j}}$ must go to infinity along a line parallel to the $t$-axis as $j$ $%
\rightarrow $ $\infty .$On the other hand, we compute%
\begin{eqnarray}
0 &=&\Theta
(v_{j})=dt(v_{j})+%
\sum_{k=1}^{n}(x_{k}(X_{L_{j}})dy_{k}-y_{k}(X_{L_{j}})dx_{k})(v_{j})
\label{vec-2} \\
&\rightarrow &\pm \infty +\text{ a bounded number}  \notag
\end{eqnarray}

\noindent as $j$ $\rightarrow $ $\infty $ since $p_{j}$ goes to a (finite)
point $p_{\infty }$ and the projection of $X_{L_{j}}$ on the $xy$-hyperplane
is a fixed vector. The contradiction obtained from (\ref{vec-2}) concludes
(a).

\endproof%

\bigskip

\section{An ODE system}

From (\ref{intcon}) (Proposition 4.2 in \cite{CCHY}), we have the following
equations for $(\alpha ,$ $k):$%
\begin{eqnarray}
e_{n}k &=&(l-2k)\alpha ,  \label{2-1} \\
e_{n}\alpha &=&k^{2}-\alpha ^{2}-kl  \notag
\end{eqnarray}

\noindent on an umbilic hypersurface $\Sigma $ of $H_{n}.$ Let $s$ be the
unit speed parameter for integral curves of $e_{n}.$ Write $k^{\prime }$ $=$ 
$\frac{dk}{ds}$ $=$ $e_{n}k$ and $\alpha ^{\prime }$ $=$ $\frac{d\alpha }{ds}
$ $=$ $e_{n}\alpha .$ Then we can write (\ref{2-1}) as follows:%
\begin{eqnarray}
k^{\prime } &=&(l-2k)\alpha ,  \label{2-1-0} \\
\alpha ^{\prime } &=&k^{2}-\alpha ^{2}-kl.  \notag
\end{eqnarray}%
\noindent Suppose $\sigma _{i,n}$ $=$ $c,$ a constant. By (\ref{1.1}), this
means that $k$ and $l$ satisfy the following relation: 
\begin{equation}
\binom{2n-2}{i-1}lk^{i-1}+\binom{2n-2}{i}k^{i}=c.  \label{2-1-1}
\end{equation}

\noindent Consider the points where $l-2k$ $=$ $0.$ We then get 
\begin{equation}
k^{i}=\frac{i}{2n+i-1}\frac{c}{\binom{2n-2}{i-1}}  \label{2-2}
\end{equation}

\noindent for $1$ $\leq $ $i$ $\leq $ $2n-1.$ For $c$ $>$ $0$ we have one
real root $k_{c,2}$ to (\ref{2-2}) if $i$ is odd while two real roots $\pm
k_{c,2}$ to (\ref{2-2}) if $i$ is even, where%
\begin{equation*}
k_{c,2}=[\frac{i}{2n+i-1}\frac{c}{\binom{2n-2}{i-1}}]^{1/i}.
\end{equation*}

\noindent We claim that the straight line $k$ $=$ $k_{c,2}$ (two straight
lines $k$ $=$ $\pm k_{c,2},$ resp.) in the $(\alpha ,$ $k)$-plane is a
solution (are solutions, resp.) to (\ref{2-1-0}) for $i$ odd (even, resp.).
When $k$ $=$ $k_{c,2}$ ($k$ $=$ $\pm k_{c,2},$ resp.)$,$ we have $l-2k$ $=$ $%
0$ and hence $k$ $=$ $k_{c,2}$ ($k$ $=$ $\pm k_{c,2},$ resp.) satisfies the
equation $k^{\prime }$ $=$ $(l-2k)\alpha $ $=$ $0.$ We can solve the second
equation of (\ref{2-1-0}) for $\alpha :$%
\begin{eqnarray*}
\alpha (s) &=&k_{c,2}\tan \{\tan ^{-1}(\frac{\alpha (s_{0})}{k_{c,2}}%
)+k_{c,2}(s_{0}-s)\} \\
(\alpha (s) &=&\pm k_{c,2}\tan \{\tan ^{-1}(\frac{\alpha (s_{0})}{\pm k_{c,2}%
})\pm k_{c,2}(s_{0}-s)\},\text{ resp.)}
\end{eqnarray*}

\noindent where $s_{0}$ $<$ $s.$ For $i$ odd (even, resp.), in $%
R^{2}\backslash \{k$ $=$ $k_{c}\}$ ($R^{2}\backslash \{k$ $=$ $\pm
k_{c,2}\}, $ resp.)$,$ there are two (three, resp.) stationary points on the 
$k$-axis ($\alpha $ $=$ $0,$ $k(k-l)=0$) for the ODE system (\ref{2-1}):%
\begin{eqnarray*}
&&(0,0)\text{ and }(0,k_{c,1}) \\
&&((0,0)\text{ and }(0,\pm k_{c,1}),\text{ resp.)}
\end{eqnarray*}

\noindent where $k_{c,1}$ is a positive number such that $l-k$ $=$ $0.$ By (%
\ref{2-1-1}) we get%
\begin{equation*}
k^{i}=\frac{i}{2n-1}\frac{c}{\binom{2n-2}{i-1}},
\end{equation*}

\noindent and hence

\begin{equation*}
k_{c,1}=[\frac{i}{2n-1}\frac{c}{\binom{2n-2}{i-1}}]^{1/i}.
\end{equation*}

We are going to prove a similar result as in Lemma 6.1 of \cite{CCHY}.

\bigskip

\textbf{Theorem 3.1}. \textit{Assume }$\sigma _{i,n}$\textit{\ }$=$\textit{\ 
}$c,$\textit{\ a positive constant.} \textit{For }$i$\textit{\ odd and }$i$ $%
\geq $ $3$ \textit{(}$i$ \textit{even and }$i$ $\geq $ $2$\textit{, resp.),
given any initial point }$p_{0}$\textit{\ }$=$\textit{\ }$(\alpha
_{0},k_{0}) $\textit{\ }

\textit{(a) Suppose} $p_{0}$ $\in $\textit{\ }$\{k>k_{c,2}\}\backslash
\{(0,k_{c,1})\}$\textit{\ (}$\{|k|>k_{c,2}\}\backslash \{(0,\pm k_{c,1})\},$%
\textit{\ resp.), there passes a unique periodic orbit }$\gamma $ $\subset $%
\textit{\ }$\{k>k_{c,2}\}\backslash \{(0,k_{c,1})\}$\textit{\ (}$%
\{|k|>k_{c,2}\}\backslash \{(0,\pm k_{c,1})\},$\textit{\ resp.), described
by }$(\alpha (s),k(s)),$\textit{\ }$0$\textit{\ }$\leq $\textit{\ }$s$%
\textit{\ }$\leq $\textit{\ }$s_{0},$\textit{\ which is a solution to the
ODE system (\ref{2-1-0}), with }$\alpha (s_{0})$\textit{\ }$=$\textit{\ }$%
\alpha (0)$\textit{\ }$=$\textit{\ }$\alpha _{0}$\textit{\ and }$k(s_{0})$%
\textit{\ }$=$\textit{\ }$k(0)$\textit{\ }$=$\textit{\ }$k_{0}.$\textit{\ }

\textit{(b) Suppose} $p_{0}$ $\in $\textit{\ }$\{0<k<k_{c,2}\}$ \textit{(}$%
\{|k|<k_{c,2},k\neq 0\}\backslash \{(0,0)\},$\textit{\ resp.),} \textit{%
there passes a unique arc-type orbit }$\gamma $ $\subset $\textit{\ }$%
\{0<k<k_{c,2}\}$ \textit{(}$\{|k|<k_{c,2}\}\backslash \{(0,0)\},$\textit{\
resp.),} \textit{described by }$(\alpha (s),k(s)),$\textit{\ for }$s$\textit{%
\ }$\in $ $(-s_{1},+s_{2})$ \textit{with} $s_{1}$ $>$ $0,$ $s_{2}$ $>$ $0,$%
\textit{\ which is a solution to the ODE system (\ref{2-1-0}), with }$%
(\alpha (0),k(0))$ $=$ $p_{0}$ \textit{and }$\lim_{s\rightarrow
+s_{2}}(\alpha (s),k(s))$\textit{\ }$=$\textit{\ }$(-\alpha _{1},0)$\textit{%
\ and }$\lim_{s\rightarrow -s_{1}}(\alpha (s),k(s))$\textit{\ }$=$\textit{\ }%
$(+\alpha _{1},0)$ \textit{for some }$\alpha _{1}$ $>$ $0.$ \textit{%
Moreover, we have}%
\begin{equation*}
\lim_{s\rightarrow -s_{1}}\frac{dk}{ds}=+\infty ,\lim_{s\rightarrow +s_{2}}%
\frac{dk}{ds}=-\infty .
\end{equation*}

\textit{(c) Suppose} $p_{0}$ $\in $\textit{\ }$\{k<0\}$ \textit{there passes
a unique arc-type orbit }$\gamma $ $\subset $\textit{\ }$\{k<0\},$ \textit{%
described by }$(\alpha (s),k(s)),$\textit{\ for all }$s$\textit{\ }$\in $ $%
R, $\textit{\ which is a solution to the ODE system (\ref{2-1-0}), with }$%
(\alpha (0),k(0))$ $=$ $p_{0}$ \textit{and }$\lim_{s\rightarrow \pm \infty
}(\alpha (s),k(s))$\textit{\ }$=$\textit{\ }$(\pm \alpha _{2},0)$\textit{\
for some} $\alpha _{2}$ $>$ $0.$

\textit{Moreover, }$\gamma $\textit{\ is symmetric with respect to the }$k$%
\textit{-axis, i.e., }$(\alpha ,k)$\textit{\ }$\in $\textit{\ }$\gamma $%
\textit{\ implies }$(-\alpha ,k)$\textit{\ }$\in $\textit{\ }$\gamma .$

\textit{\bigskip }

There are important curves $\Psi _{0}^{\pm }$ in $(\alpha ,k)$-plane,
defined by $k^{2}-\alpha ^{2}-kl$ $=$ $0$ (where $\alpha ^{\prime }$ $=$ $0$%
), in $k\gtrless $ 0 regions, resp. (see Figure 3.1). Recall that $l$ is a
function of $k$ by (\ref{2-1-1}):%
\begin{equation}
l=\frac{c}{\binom{2n-2}{i-1}k^{i-1}}-\frac{2n-i-1}{i}k  \label{3-0}
\end{equation}

\noindent for $1\leq i\leq 2n-1$ (we have used the convention $\binom{m}{0}$ 
$=$ $1).$ Compute%
\begin{equation}
\frac{dl}{dk}=-(i-1)\frac{c}{\binom{2n-2}{i-1}k^{i}}-\frac{2n-i-1}{i}.
\label{3-1}
\end{equation}

\noindent Expressing $\alpha $ in $k$ from $k^{2}-\alpha ^{2}-kl$ $=$ $0,$
we get%
\begin{equation*}
\alpha =k\sqrt{\frac{2n-1}{i}-\frac{c}{\binom{2n-2}{i-1}k^{i}}}
\end{equation*}

\noindent in the region: $\alpha $ $>$ $0$ and $k$ $>$ $k_{c,1}.$ Now write $%
a$ $=$ $\frac{2n-1}{i},$ $b$ $=$ $\frac{c}{\binom{2n-2}{i-1}}.$ A direct
computation shows%
\begin{eqnarray*}
\frac{d^{2}\alpha }{dk^{2}} &=&\frac{bik^{-i-1}}{2\sqrt{a-bk^{-i}}}+\frac{%
-bi^{2}k^{-i-1}}{2\sqrt{a-bk^{-i}}}-\frac{b^{2}i^{2}k^{-2i-1}}{%
4(a-bk^{-i})^{3/2}} \\
&<&0\text{ \ \ (since }i^{2}\geq i\text{).}
\end{eqnarray*}

\noindent From the above argument and symmetry in $k$-axis, we have shown
the following result.

\bigskip

\textbf{Lemma 3.1}. $\Psi _{0}^{+}$\textit{\ is convex upward in the region: 
}$\alpha $\textit{\ }$>$\textit{\ }$0$\textit{\ (}$\alpha <0,$\textit{\
resp.)} \textit{and }$k$\textit{\ }$>$\textit{\ }$k_{c,1},$\textit{\ viewed
as a graph over }$k$\textit{-axis. Moreover, }$\alpha $ \textit{goes like} $%
\sqrt{\frac{2n-1}{i}}k$ \textit{(}$-\sqrt{\frac{2n-1}{i}}k$\textit{\ in the }%
$\alpha <0$\textit{\ region, resp.) as }$k$\textit{\ tends to +}$\infty .$

\bigskip

We are going to compute curvature of a solution curve $(\alpha (s),k(s))$ to
(\ref{2-1-0}):%
\begin{equation*}
\frac{d}{ds}(\tan ^{-1}\frac{k^{\prime }}{\alpha ^{\prime }})=\frac{%
k^{\prime \prime }\alpha ^{\prime }-\alpha ^{\prime \prime }k^{\prime }}{%
(\alpha ^{\prime })^{2}+(k^{\prime })^{2}}.
\end{equation*}

Taking derivatives of (\ref{2-1-0}) and substituting into $k^{\prime \prime
}\alpha ^{\prime }-\alpha ^{\prime \prime }k^{\prime },$ we get%
\begin{eqnarray}
&&k^{\prime \prime }\alpha ^{\prime }-\alpha ^{\prime \prime }k^{\prime }
\label{3.2} \\
&=&-(2k-l)\{k^{2}\alpha ^{2}[(\frac{2n-1}{i^{2}})(2n+3i-1)  \notag \\
&&+k^{-i}\tilde{\sigma}_{i,n}(\frac{-4n+i^{2}-3i+2}{i})+k^{-2i}(\tilde{\sigma%
}_{i,n})^{2}]  \notag \\
&&+\alpha ^{4}[\frac{2n-i-1}{i}+(i-1)\tilde{\sigma}_{i,n}k^{-i}]+(\alpha
^{\prime })^{2}\}.  \notag
\end{eqnarray}

\noindent In deriving (\ref{3.2}) we have used (\ref{3-1}) and let%
\begin{equation*}
\tilde{\sigma}_{i,n}=\frac{\sigma _{i,n}}{\binom{2n-2}{i-1}}
\end{equation*}

\noindent (see the expression of $\sigma _{i,n}$ in (\ref{1.1})). Denote the
main term in (\ref{3.2}) by%
\begin{eqnarray}
\Pi &:&=k^{2}\alpha ^{2}[(\frac{2n-1}{i^{2}})(2n+3i-1)  \label{3.3} \\
&&+k^{-i}\tilde{\sigma}_{i,n}(\frac{-4n+i^{2}-3i+2}{i})+k^{-2i}(\tilde{\sigma%
}_{i,n})^{2}]  \notag \\
&&+\alpha ^{4}[\frac{2n-i-1}{i}+(i-1)\tilde{\sigma}_{i,n}k^{-i}]+(\alpha
^{\prime })^{2}  \notag
\end{eqnarray}

\noindent where $\alpha ^{\prime }$ in the last term should mean $%
k^{2}-\alpha ^{2}-kl$ by (\ref{2-1-0}). So from (\ref{3.2}) and (\ref{3.3}),
we can write%
\begin{eqnarray}
\frac{d^{2}\alpha }{dk^{2}} &=&\frac{d}{dk}(\frac{\alpha ^{\prime }}{%
k^{\prime }})  \label{3.4} \\
&=&\frac{\alpha ^{\prime \prime }k^{\prime }-\alpha ^{\prime }k^{\prime
\prime }}{(k^{\prime })^{3}}=\frac{(2k-l)\Pi }{(k^{\prime })^{3}}  \notag \\
&=&-\frac{\Pi }{(2k-l)^{2}\alpha ^{3}}  \notag
\end{eqnarray}

\noindent by (\ref{2-1-0}).

\bigskip

\proof
\textbf{(of Theorem 3.1)} If the initial point $p_{0}$ is in $k$-axis with $%
k_{0}$ $>$ $k_{c,1},$ then $k^{\prime }(0)$ $=$ $0$ and $\alpha ^{\prime
}(0) $ $>$ $0.$ So according to the ODE (\ref{2-1-0}), $p_{0}$ is moving
into the region surrounded by the line segment $\{0\}x[k_{c,1},\infty )$ and
the curve $\Psi _{0}^{+}$ $\cap $ $\{\alpha >0\}$ in small positive time.
Call this region $I^{+}$ (see Figure 3.1)$.$ Since the solution to (\ref%
{2-1-0}) is symmetric with respect to the $k$-axis, we need only to discuss
the $\alpha $-positive part. We may assume $i$ $\geq $ $2$ in the following
argument (for the case $i$ $=$ $1,$ $\sigma _{1,n}$ is the $p$-mean
curvature. We refer the reader to \cite{CCHY}).

\begin{figure}[h]
\includegraphics[width= 10.4cm]{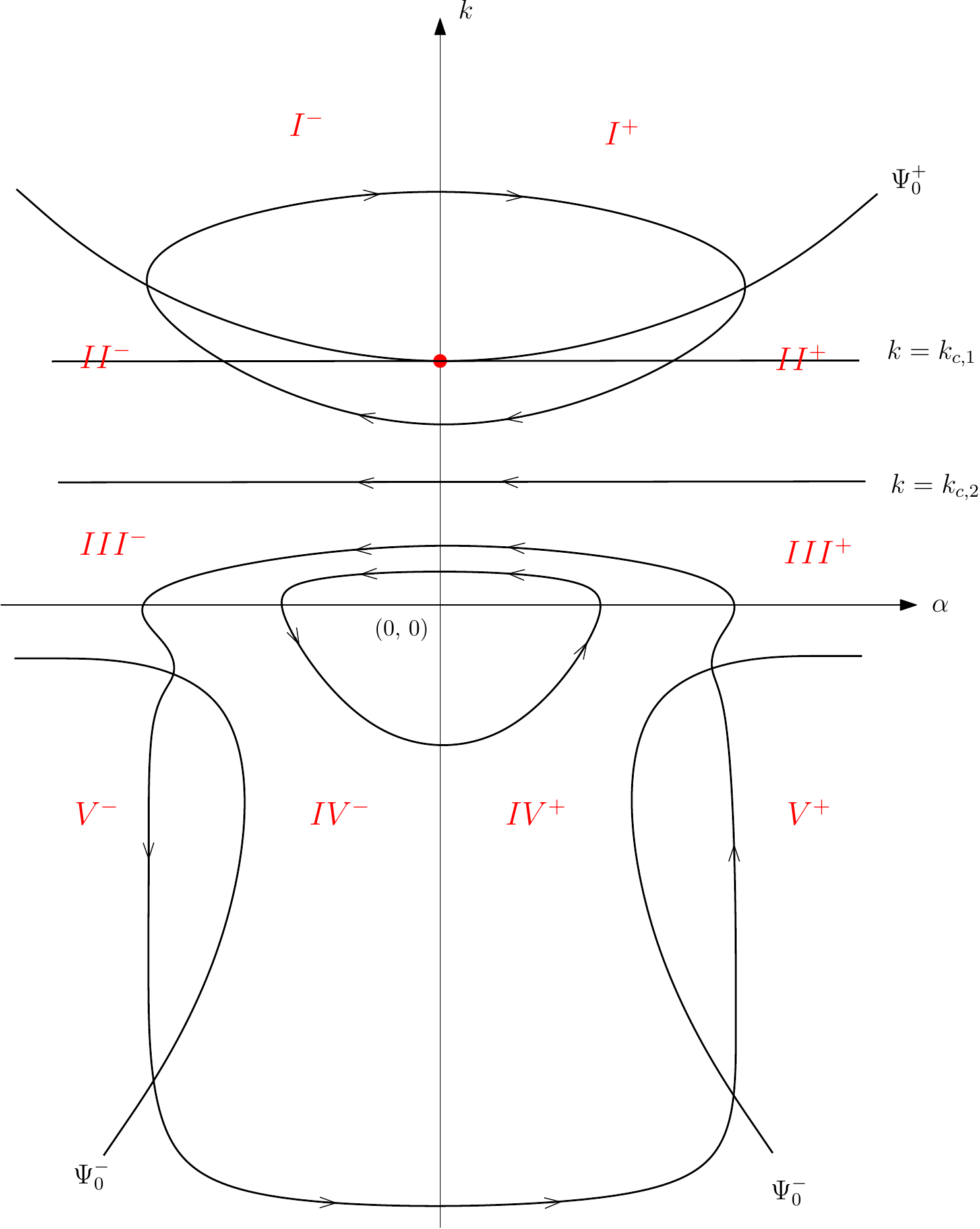} 
\caption{Figure 3.1: $c>0$, odd $i \ge 3$}
\end{figure}

\textbf{Case 1}. Suppose $p_{0}$ $\in $ Region $I^{+}.$ From the equation (%
\ref{2-1-0}) we learn that $k^{\prime }$ $=$ $(l-2k)\alpha $ $<$ $0$ since $%
\alpha $ $>$ $0$ and $l-2k$ $<$ $0$ for $k$ $>$ $k_{c,1}.$ On the other
hand, $\alpha ^{\prime }$ $=$ $k^{2}-\alpha ^{2}-kl$ $>$ $0$ since $p_{0}$
is above the curve $\Psi _{0}^{+}.$ Altogether we have%
\begin{equation}
\frac{dk}{d\alpha }=\frac{k^{\prime }}{\alpha ^{\prime }}=\frac{(l-2k)\alpha 
}{k^{2}-\alpha ^{2}-kl}<0  \label{2-3}
\end{equation}

\noindent in Region $I^{+}.$ It follows from (\ref{2-3}) and Lemma 3.1 that
the solution $(\alpha (s),k(s))$ to (\ref{2-1-0}) with \textit{\ }$(\alpha
(0),k(0))$ $=$ $p_{0}$ must hit $\Psi _{0}^{+}$ in finite positive time.
That is, for some $s_{1}$ $>$ $0,$ $(\alpha (s_{1}),k(s_{1}))$ $\in $ $\Psi
_{0}^{+}.$ Next what happens as $s\rightarrow -\infty ?$ In view of (\ref%
{2-3}) we have $(\alpha (s),k(s))$ $\in $ $\Omega _{0}$ $:=$ [0,$\alpha
_{0}]\times \lbrack k_{0},\infty )$ for $s$ $\leq $ $0.$ There exists large $%
k_{1}$ $>>$ $k_{c,1}$ such that for $k$ $>$ $k_{1},$ $|\alpha |$ $\leq $ $%
\alpha _{0},$ we have $2k-l$ $=$ $O(k)$, $k^{2}-kl$ $\geq $ $k^{2}$, and
hence%
\begin{eqnarray}
|\frac{dk}{d\alpha }| &=&|\frac{k^{\prime }}{\alpha ^{\prime }}|  \label{2-4}
\\
&=&|\frac{(l-2k)\alpha }{k^{2}-\alpha ^{2}-kl}|  \notag \\
&\leq &\frac{C_{1}}{k}\leq C_{2}  \notag
\end{eqnarray}

\noindent in $\Omega _{0}.$ Therefore $(\alpha (s),k(s))$ must hit \{0\}$%
\times (k_{c,1},\infty )$ in negative finite time. That is, for some $s_{2}$ 
$<$ $0,$ $(\alpha (s_{2}),k(s_{2}))$ $\in $ \{0\}$\times (k_{c,1},\infty ).$

If $p_{0}$ $\in $ $\Psi _{0}^{+}$ $\cap $ $\{\alpha >0\},$ $\alpha ^{\prime
}(0)$ $=$ $0$ and $k^{\prime }(0)$ $=$ $(l-2k)(0)\alpha (0)$ $<$ $0$ since $%
(l-2k)(0)$ $<$ $0,$ $\alpha (0)$ $>$ $0.$ So by Lemma 3.1, $p_{0}$ will move
into the region $II^{+}$ surrounded by the curve $\Psi _{0}^{+}$ $\cap $ $%
\{\alpha >0\},$ line segment \{0\}$\times \lbrack k_{c,2},k_{c,1}],$ and the
horizontal line segment $L_{c,2}^{+}$ $:=$ $\{(\alpha ,k_{c,2})$ $:$ $\alpha 
$ $\geq $ $0\}$ in small positive time.

\textbf{Case 2}. Suppose $p_{0}$ $\in $ Region $II^{+}.$ Write Region $%
II^{+} $ $=$ Region $II_{1}$ $\cup $ the horizontal line segment $%
L_{c,1}^{+} $ $:=$ $\{(\alpha ,k_{c,1})$ $:$ $\alpha $ $\geq $ $0\}$ $\cup $
Region $II_{2}.$ Here Region $II_{1}$ denotes the region surrounded by the
curve $\Psi _{0}^{+}$ $\cap $ $\{\alpha >0\}$ and $L_{c,1}^{+}$ while Region 
$II_{2} $ denotes the region surrounded by $L_{c,1}^{+},$ line segment \{0\}$%
\times \lbrack k_{c,2},k_{c,1}],$ and $L_{c,2}^{+}.$

Observe that the "discriminant" 
\begin{equation*}
4(\frac{2n-1}{i})(2n+3i-1)>(\frac{-4n+i^{2}-3i+2}{i})^{2}
\end{equation*}

\noindent holds for $1$ $\leq $ $i$ $\leq $ $2n-1,$ $n$ $\geq $ $2$ by
elementary algebra. From (\ref{3.3}) we then have the following estimate:%
\begin{eqnarray}
&&\frac{\Pi }{(2k-l)^{2}\alpha ^{3}}  \label{3.5} \\
&\geq &\frac{C_{1}(k^{2}\alpha ^{2}+\alpha ^{4})}{(2k-l)^{2}\alpha ^{3}}=%
\frac{C_{1}(k^{2}+\alpha ^{2})}{(2k-l)^{2}\alpha }  \notag \\
&\geq &C_{2}\frac{k\alpha }{k^{2}\alpha }=C_{2}\frac{1}{k}  \notag
\end{eqnarray}

\noindent for $k$ $\geq $ $k_{0}$ $>$ $k_{c,2}$ and $\alpha $ $>$ $0.$ Here
we have used 
\begin{equation}
l-2k=\frac{c}{\binom{2n-2}{i-1}}k^{1-i}-\frac{2n+i-1}{i}k  \label{3.5.1}
\end{equation}%
\noindent deduced from (\ref{3-0}). The positive constants $C_{1},$ $C_{2}$
depend on $n,i,$ $c,$ and $k_{0}.$ From (\ref{3.4}) and (\ref{3.5}), we have%
\begin{eqnarray}
\frac{d\alpha }{dk}(k)-\frac{d\alpha }{dk}(k_{1}) &=&\int_{k_{1}}^{k}\frac{%
d^{2}\alpha }{dk^{2}}dk  \label{3.6} \\
&\leq &-C_{2}\int_{k_{1}}^{k}\frac{1}{k}dk  \notag \\
&=&-C_{2}\log \frac{k}{k_{1}}  \notag
\end{eqnarray}

\noindent for $k$ $>$ $k_{1}$ $\geq $ $k_{0}$ ($>$ $k_{c,2}$ $>$ $0).$ The
solution $(\alpha (s),k(s))$ to (\ref{2-1-0}) with \textit{\ }$(\alpha
(0),k(0))$ $=$ $p_{0}$ $\in $ $II^{+}$ must hit $\Psi _{0}^{+}$ $\cap $ $%
\{\alpha >0\}$ in finite negative time. Otherwise as $s$ $\rightarrow $ $%
-\infty ,$ $(\alpha (s),k(s))$ stays in $II^{+}$ while $k$ goes to $+\infty $
or some finite $k_{2}.$ In the first case, we let $k$ $\rightarrow $ +$%
\infty $ in (\ref{3.6}) to obtain%
\begin{equation*}
0-\frac{d\alpha }{dk}(k_{1})\leq \lim_{k\rightarrow \infty }\inf \frac{%
d\alpha }{dk}(k)-\frac{d\alpha }{dk}(k_{1})\leq -\infty ,
\end{equation*}

\noindent a contradiction. Here we have used the fact that $\frac{d\alpha }{%
dk}$ $>$ $0$ in $II^{+}.$ In the second case, we have $\lim_{k\rightarrow
k_{2}-}\frac{d\alpha }{dk}$ $=$ $+\infty $ while, from (\ref{3.6}), there
holds 
\begin{equation*}
+\infty =\lim_{k\rightarrow k_{2}-}\frac{d\alpha }{dk}(k)-\frac{d\alpha }{dk}%
(k_{1})\leq -C_{2}\log \frac{k_{2}}{k_{1}},
\end{equation*}

\noindent a contradiction. On the other hand, it is not possible that as $s$ 
$\rightarrow $ $+\infty ,$ the solution $(\alpha (s),k(s))$ stays in Region $%
II_{1}.$ The reason is that $\lim_{k\rightarrow k_{c,1}+}\frac{d\alpha }{dk}%
(k)$ $=$ $+\infty $ in this case while, from (\ref{3.6}) (taking $k_{1}$
such that $k_{c,1}$ $>$ $k_{1}$ $>$ $k_{c,2})$, there holds%
\begin{equation*}
+\infty =\lim_{k\rightarrow k_{c,1}+}\frac{d\alpha }{dk}(k)-\frac{d\alpha }{%
dk}(k_{1})\leq -C_{2}\log \frac{k_{c,1}}{k_{1}},
\end{equation*}

\noindent a contradiction. Observe that both $\alpha ^{\prime }(s)$ and $%
k^{\prime }(s)$ are negative for $(\alpha (s),k(s))$ $\in $ $II^{+}.$ The
solution $(\alpha (s),k(s))$ must hit line segment \{0\}$\times
(k_{c,2},k_{c,1})$ in finite positive time since it cannot hit $L_{c,2}^{+}$
(which is a solution orbit) by the uniqueness of solutions to an ODE.

\textbf{Case 3. }Let us look at Region $III^{+},$ the region surrounded by $%
L_{c,2}^{+},$ line segment \{0\}$\times \lbrack 0,k_{c,2}],$ and positive $%
\alpha $-axis. Observe that in $III^{+},$ the main term (see (\ref{3.3})) $%
\Pi $ $>$ $0$ since $c$ $>$ $0,$ $k$ $>$ $0.$ Also $\alpha $ $>$ $0$ in $%
III^{+},$ and hence we have%
\begin{equation*}
\frac{d^{2}\alpha }{dk^{2}}=-\frac{\Pi }{(2k-l)^{2}\alpha ^{3}}<0.
\end{equation*}

\noindent It follows that the graph $\alpha $ $=$ $\alpha (k)$ is convex
upward, so any solution $(\alpha (s),k(s))$ to (\ref{2-1-0}) in $III^{+}$
must hit line segment \{0\}$\times \lbrack 0,k_{c,2}]$ in finite positive
time while approaching to $(\alpha _{0},0)$ as $s$ $\rightarrow $ $-s_{1}$
for (finite) $s_{1}$ $>$ $0.$ The reason is that it cannot hit $L_{c,2}^{+}$
by the uniqueness of solutions to an ODE since $L_{c,2}^{+}$ is a solution
orbit. On the other hand, we have%
\begin{eqnarray}
\frac{dk}{ds} &=&(l-2k)(s)\alpha (s)  \label{3.6.1} \\
&=&[\frac{c}{\binom{2n-2}{i-1}}k^{1-i}(s)-\frac{2n+i-1}{i}k(s)]\alpha (s) 
\notag \\
&&\text{approximates }\frac{c}{\binom{2n-2}{i-1}}k^{1-i}\alpha _{0}  \notag
\end{eqnarray}

\noindent for $(\alpha ,k)$ near $(\alpha _{0},0)$ and $i$ $\geq $ $2$. It
follows that $k^{i-1}dk$ approximates $\tilde{c}\alpha _{0}ds$ with $\tilde{c%
}\alpha _{0}$ $>$ $0,$ and hence after integration, we learn that $k(s)$
goes to zero only when $s$ goes to a finite number $-s_{1}$. By (\ref{3.6.1}%
) we get%
\begin{equation*}
\lim_{s\rightarrow -s_{1}}\frac{dk}{ds}=+\infty .
\end{equation*}

Next, we will discuss regions $IV^{+}$ and $V^{+}$ for odd $i$ $\geq $ $3$
(see Figure 3.1). For even $i$ $\geq $ $2,$ we are done since the solution
to (\ref{2-1-0}) is also symmetric with respect to the $\alpha $-axis (see
Figure 3.2).

\begin{figure}[h]
\includegraphics[width= 10.4cm]{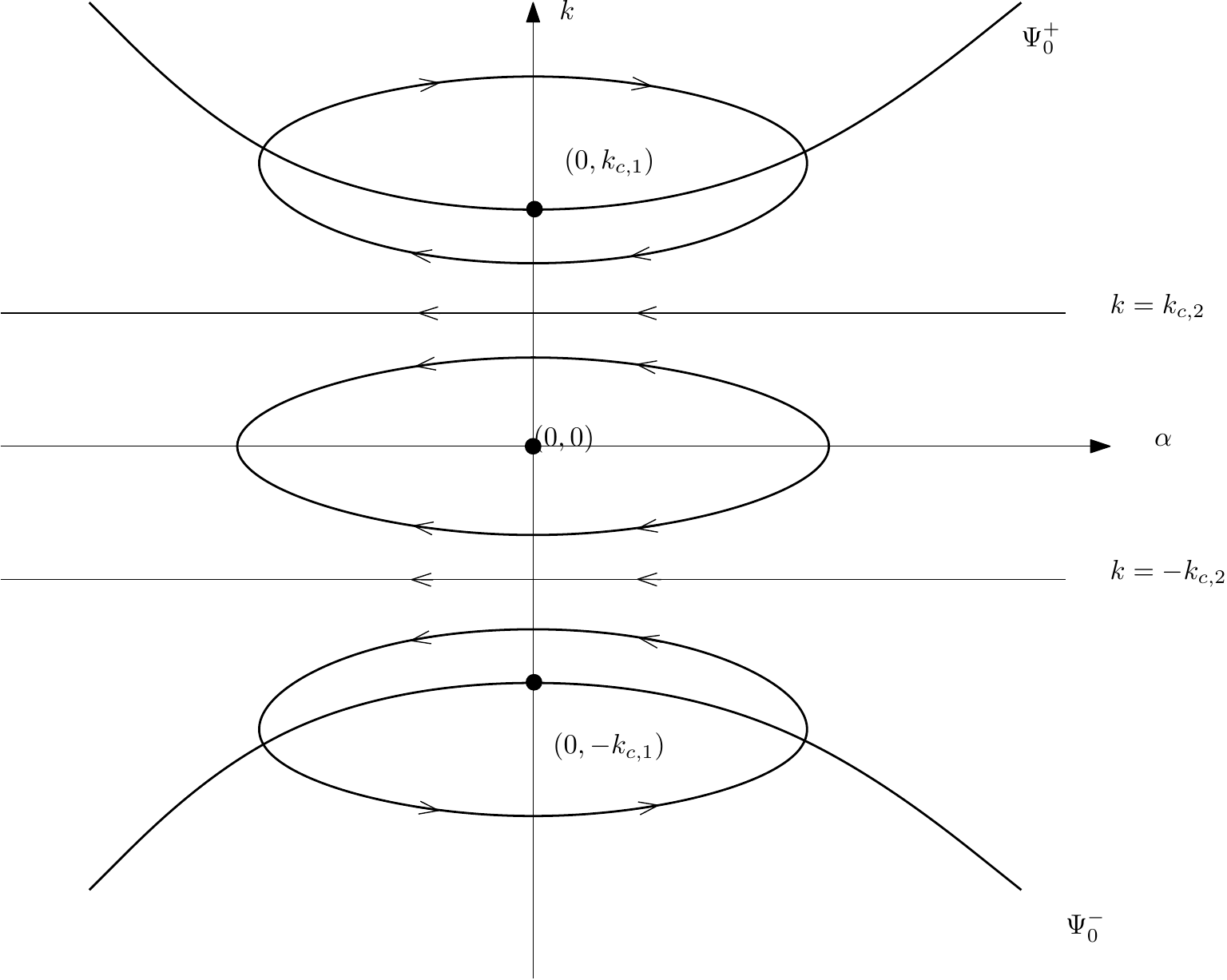} 
\caption{Figure 3.2: $c>0$, even $i\geq 2$}
\end{figure}

\textbf{Case 4}. Look at the solutions to (\ref{2-1-0}) in region $V^{+}$:
surrounded by $\Psi _{0}^{-}$ $\cap $ $\{\alpha >0\},$ in which $\alpha
^{\prime }$ $<$ $0$ (see Figure 3.1)$.$ We estimate the main term $\Pi $
(see (\ref{3.3})): there is a large negative $k_{0}$ such that for $k$ $\leq 
$ $k_{0}$ $<$ $0$ there holds%
\begin{equation}
\Pi \geq C_{1}(k^{2}\alpha ^{2}+\alpha ^{4})  \label{3.7}
\end{equation}

\noindent where the positive constant $C_{1}$ depends on $c,$ $k_{0},$ $n,$
and $i.$ From (\ref{3.4}), (\ref{3.7}), and (\ref{3.5.1}), we then have%
\begin{eqnarray}
\frac{d^{2}\alpha }{dk^{2}} &=&-\frac{\Pi }{(2k-l)^{2}\alpha ^{3}}
\label{3.8} \\
&\leq &-\frac{C_{1}(k^{2}\alpha ^{2}+\alpha ^{4})}{(2k-l)^{2}\alpha ^{3}} 
\notag \\
&\leq &-C_{1}\frac{k^{2}+\alpha ^{2}}{(2k-l)^{2}\alpha }  \notag \\
&\leq &-C_{2}\frac{|k|\alpha }{k^{2}\alpha }=-C_{2}\frac{1}{|k|}  \notag
\end{eqnarray}

\noindent for $k$ $\leq $ $k_{0}$ $<$ $0$ (may need to take larger $k_{0}).$
Here positive constants $C_{1},$ $C_{2}$ depend only on $c,$ $k_{0},$ $n,$
and $i.$ Integrating (\ref{3.8}) gives%
\begin{equation}
\frac{d\alpha }{dk}(k_{0})-\frac{d\alpha }{dk}(k)\leq C_{2}\log |\frac{k_{0}%
}{k}|.  \label{3.9}
\end{equation}

\noindent We claim the solution $(\alpha (s),k(s))$ to (\ref{2-1-0}) with 
\textit{\ }$(\alpha (0),k(0))$ $=$ $p_{0}$ $\in $ $V^{+}$ with $k(0)$ $\leq $
$k_{0}$ must hit the curve $\Psi _{0}^{-}$ $\cap $ $\{\alpha >0\}$ in finite
negative time. Otherwise we have either $\lim_{k\rightarrow -\infty }\frac{%
d\alpha }{dk}(k)$ $\leq $ $0$ or $\lim_{k\rightarrow k_{1}+}\frac{d\alpha }{%
dk}(k)$ $=$ $-\infty $ since $\frac{d\alpha }{dk}$ is monotonically
decreasing as $k$ decreases (as $s$ goes to $-\infty )$. But from (\ref{3.9}%
) we get%
\begin{equation}
\frac{d\alpha }{dk}(k_{0})\leq \frac{d\alpha }{dk}(k_{0})-\lim_{k\rightarrow
-\infty }\frac{d\alpha }{dk}(k)\leq C_{2}(-\infty )  \label{3.9.1}
\end{equation}

\noindent or 
\begin{equation}
-\infty =\lim_{k\rightarrow k_{1}+}\frac{d\alpha }{dk}(k)-\frac{d\alpha }{dk}%
(k_{0})\leq C_{2}\log |\frac{k_{0}}{k_{1}}|.  \label{3.9.2}
\end{equation}

\noindent Both (\ref{3.9.1}) and (\ref{3.9.2}) give contradiction. Next
suppose $k(0)$ $>$ $k_{0}.$ If $k(s_{1})$ $\leq $ $k_{0}$ for $s_{1}$ $<$ $%
0, $ we run the above argument again to show the solution $(\alpha (s),k(s))$
must hit the curve $\Psi _{0}^{-}$ $\cap $ $\{\alpha >0\}$ in finite
negative time. In case $k(s)$ $>$ $k_{0}$ for all $s$ $<$ $0,$ we must have%
\begin{equation}
\lim_{k\rightarrow k_{2}+}\frac{d\alpha }{dk}(k)=-\infty \text{ and }%
\lim_{k\rightarrow k_{2}+}\alpha (k)=+\infty  \label{3.9.3}
\end{equation}

\noindent for some $k_{2}$ $\geq $ $k_{0}.$ There exists a large $\alpha
_{0} $ such that for $\alpha $ $\geq $ $\alpha _{0}$, $k$ close to $k_{2},$
we have the following estimate:%
\begin{eqnarray}
\frac{d^{2}\alpha }{dk^{2}} &=&-\frac{\Pi }{(2k-l)^{2}\alpha ^{3}}
\label{3.10} \\
&\leq &-\frac{k^{-i}\frac{c}{\binom{2n-2}{i-1}}[k^{2}\alpha ^{2}(\frac{%
i^{2}+2}{i})+\alpha ^{4}(i-1)]}{(2k-l)^{2}\alpha ^{3}}  \notag \\
&\leq &C_{3}\alpha  \notag
\end{eqnarray}

\noindent for some poitive constant $C_{3}$ depending on $k_{2},$ $\alpha
_{0},$ $c,$ $n,$ $i.$ Note that $k^{-i}$ $<$ $0$ since $k$ $<$ $0$ and $i$
is odd. For the first inequality in (\ref{3.10}), we have dropped positive
terms in $\Pi .$ Multiplying (\ref{3.10}) by $2\frac{d\alpha }{dk}$ and
integrating from $k$ to $k_{3},$ $k_{3}$ $>$ $k$ $>$ $k_{2},$ we obtain%
\begin{equation*}
(\frac{d\alpha }{dk})^{2}(k_{3})-(\frac{d\alpha }{dk})^{2}(k)\geq
C_{3}(\alpha ^{2}(k_{1})-\alpha ^{2}(k)).
\end{equation*}

\noindent Hence%
\begin{equation}
(\frac{d\alpha }{dk})^{2}(k)\leq C_{3}\alpha ^{2}(k)+C_{4}  \label{3.11}
\end{equation}

\noindent where $C_{4}$ $=$ $(\frac{d\alpha }{dk})^{2}(k_{3})$ $-$ $%
C_{3}\alpha ^{2}(k_{1}).$ If $C_{4}$ $\leq $ $0,$ we have $(\frac{d\alpha }{%
dk})^{2}(k)\leq C_{3}\alpha ^{2}(k)$ from (\ref{3.11}), and hence%
\begin{eqnarray}
\log \frac{\alpha (k)}{\alpha (k_{3})} &=&-\int_{k}^{k_{3}}\frac{d\alpha }{%
\alpha }  \label{3.12} \\
&=&\int_{k}^{k_{3}}\frac{|d\alpha |}{\alpha }\leq \int_{k}^{k_{3}}\sqrt{C_{3}%
}|dk|  \notag \\
&=&\sqrt{C_{3}}(k_{3}-k)\rightarrow \sqrt{C_{3}}(k_{3}-k_{2})  \notag
\end{eqnarray}

\noindent as $k\rightarrow k_{2}.$ On the other hand, $\log \frac{\alpha (k)%
}{\alpha (k_{3})}$ on the left hand side of (\ref{3.12}) tends to +$\infty $
as $k\rightarrow k_{2}$ by (\ref{3.9.3}). We have reached +$\infty $ $\leq $
a finite number, a contradiction. In case $C_{4}$ $>$ $0,$ we can absorb it
in $\alpha ^{2}(k)$ term for $k$ large so that $(\frac{d\alpha }{dk}%
)^{2}(k)\leq C_{5}\alpha ^{2}(k)$ for $C_{5}$ $>$ $C_{3}.$ A similar
reasoning as above reaches a contradiction again. We have shown that in any
case the solution $(\alpha (s),k(s))$ must hit the curve $\Psi _{0}^{-}$ $%
\cap $ $\{\alpha >0\}$ in finite negative time.

On the other hand, we claim $(\alpha (s),k(s))$ must hit the curve $\Psi
_{0}^{-}$ $\cap $ $\{\alpha >0\}$ also in finite positive time. Observe that
as $s$ goes to $+\infty ,$ $\alpha (s)$ decreases while $k(s)$ increases. So 
$(\alpha (s),k(s))$ will hit $\Psi _{0}^{-}$ $\cap $ $\{\alpha >0\}$ in
finite positive time for topological reason.

\textbf{Case 5.} Let us look at solutions $\gamma (s)$ :$=$ $(\alpha
(s),k(s))$ to (\ref{2-1-0}), which lie in Region $IV^{+}:$ the region
surrounded by positive $\alpha $-axis, negative $k$-axis, and the curve $%
\Psi _{0}^{-}$ $\cap $ $\{\alpha >0\}$ (see Figure 3.1). Recall that $k_{0}$
is a large negative number in Case 4, such that estimate (\ref{3.8}) holds
for $k$ $\leq $ $k_{0}.$ If $k(s_{1})$ $\leq $ $k_{0}$ for some $s_{1},$
then $\gamma (s)$ must hit the negative $k$-axis in finite negative time
since $\gamma (s)$ viewed as a graph $\alpha $ $=$ $\alpha (k)$ is convex
upward by (\ref{3.8}). In the other direction, $\gamma (s)$ must hit $k$ $=$ 
$k_{0}$ line (then enters $IV^{+}$ $\cap $ $\{k$ $>$ $k_{0}\})$ or the curve 
$\Psi _{0}^{-}$ $\cap $ $\{\alpha >0\}$ (then enters region $V^{+})$ in
finite positive time.

In case $k(s_{2})$ $>$ $k_{0}$ for some $s_{2},$ $\gamma (s)$ will backward
hit negative $k$-axis or $k$ $=$ $k_{0}$ line or curve $\Psi _{0}^{-}$ $\cap 
$ $\{\alpha >0\}$ in finite negative time by topological reason (noting that 
$\alpha ^{\prime }(s)$ $>$ $0$ and $k^{\prime }(s)$ $>$ $0$ for $\gamma (s)$ 
$\in $ $IV^{+})$. In the other direction, $\gamma (s)$ will hit $\Psi
_{0}^{-}$ $\cap $ $\{\alpha >0\}$ (then enters region $V^{+})$ in finite
positive time or hit positive $\alpha $-axis in infinite positive time (the
same reason as shown in the discussion of Case 3). If both cases do not
occur, then $\gamma (s)$ goes to the infinity between positive $\alpha $%
-axis and curve $\Psi _{0}^{-}$ $\cap $ $\{\alpha >0\}.$ We claim this
cannot happen. Observe that in $IV^{+},$ we have%
\begin{eqnarray}
\alpha ^{2} &\leq &k(k-l)  \label{3.12.1} \\
&\leq &-C_{6}k^{2-i}  \notag
\end{eqnarray}
\noindent for $k$ near $0$ (noting that $i$ is odd, $i$ $\geq $ $3,$ and $k$
is negative)$.$ As $\alpha $ $\rightarrow $ $+\infty ,$ $k$ must go to $0$
by the first inequality of (\ref{3.12.1}). Besides, the situation also
forces 
\begin{equation}
\lim_{k\rightarrow 0^{-}}\frac{d\alpha }{dk}=+\infty .  \label{3.13}
\end{equation}

On the other hand, we estimate $\Pi $ $\geq $ $C_{7}\alpha ^{4}k^{-i}$ and
hence%
\begin{eqnarray}
\frac{d^{2}\alpha }{dk^{2}} &=&-\frac{\Pi }{(2k-l)^{2}\alpha ^{3}}
\label{3.14} \\
&\leq &-C_{7}\frac{\alpha ^{4}k^{-i}}{(2k-l)^{2}\alpha ^{3}}  \notag \\
&\leq &-C_{8}\frac{\alpha k^{-i}}{(k^{1-i})^{2}}=-C_{8}\frac{\alpha }{k^{2-i}%
}  \notag \\
&\leq &-C_{8}\frac{\alpha ^{2}}{k^{2-i}}\leq C_{9}  \notag
\end{eqnarray}

\noindent by (\ref{3.12.1}) for $\alpha $ large ($k$ near $0)$. Integrating (%
\ref{3.14}) gives%
\begin{equation*}
\frac{d\alpha }{dk}(k)-\frac{d\alpha }{dk}(\hat{k})\leq C_{9}(k-\hat{k})
\end{equation*}

\noindent for fixed $\hat{k}$, close to $0$ and $\hat{k}$ $<$ $k$ $<$ $0.$
Letting $k$ $\rightarrow $ $0^{-},$ we get%
\begin{equation*}
+\infty \leq C_{9}(-\hat{k})
\end{equation*}

\noindent by (\ref{3.13}), a contradiction.

The discussion for the part of $\alpha $ $<$ $0$ is similar since the
solution $\gamma (s)$ :$=$ $(\alpha (s),k(s))$ to (\ref{2-1-0}) is symmetric
with respect to the $k$-axis, i.e., $(\alpha ,k)$\ $\in $\ $\gamma $\
implies $(-\alpha ,k)$\ $\in $\ $\gamma .$

\endproof%

\bigskip

For convexity of solutions to (\ref{2-1-0}) we notice the following result
(see Figure 3.1 and Figure 3.2).

\bigskip

\textbf{Corollary 3.2}. \textit{Suppose we are in the same situation of
Theorem 3.1. For }$i$\textit{\ }$\geq $ $3$ \textit{odd, each solution to (%
\ref{2-1-0}) in the upper half plane is convex while not convex in the lower
half plane. For }$i$ $\geq $ $2$ \textit{even, each solution to (\ref{2-1-0}%
) is always convex.}

\bigskip

We now turn to the case $\sigma _{i,n}$\textit{\ }$=$\textit{\ }$c$ $=$ $0.$

\bigskip

\textbf{Theorem 3.3.} \textit{Assume }$\sigma _{i,n}$\textit{\ }$=$\textit{\ 
}$c$ $=$ $0.$ \textit{We are given any initial point }$p_{0}$\textit{\ }$=$%
\textit{\ }$(\alpha _{0},k_{0}).$

\textit{(a) Suppose }$p_{0}$\textit{\ }$=$\textit{\ }$(0,0).$\textit{\ Then }%
$(\alpha (s)$\textit{\ }$\equiv $\textit{\ }$0,$\textit{\ }$k(s)$\textit{\ }$%
\equiv $\textit{\ }$0)$\textit{\ is the unique solution to (\ref{2-1-0})
with (}$\alpha (0),$\textit{\ }$k(0))$\textit{\ }$=$\textit{\ }$(0,0).$%
\textit{\ Moreover, a connected, immersed umbilic hypersurface corresponding
to }$\alpha $\textit{\ }$\equiv $\textit{\ }$0,$\textit{\ }$k$\textit{\ }$%
\equiv $\textit{\ }$0$\textit{\ is congruent with part of the hypersurface }$%
E\times R$\textit{\ where }$E$\textit{\ is a hyperplane of }$R^{2n}.$

\textit{(b) Suppose }$p_{0}$\textit{\ }$=$\textit{\ }$(\alpha _{0},0)$%
\textit{\ with }$\alpha _{0}\neq 0.$\textit{\ Then we have the unique
solution to (\ref{2-1-0}) with (}$\alpha (0),$\textit{\ }$k(0))$\textit{\ }$%
= $\textit{\ }$(\alpha _{0},0)$\textit{\ as follows:}%
\begin{equation*}
\alpha (s)=\frac{1}{s+\alpha _{0}^{-1}},k(s)\equiv 0
\end{equation*}%
\textit{\noindent where }$s$\textit{\ }$\in $\textit{\ (}$-\alpha
_{0}^{-1},+\infty )$\textit{\ (}$s$\textit{\ }$\in $\textit{\ (}$-\infty
,-\alpha _{0}^{-1}),$\textit{\ resp.) for }$\alpha _{0}$\textit{\ }$>$%
\textit{\ }$0$\textit{\ (}$\alpha _{0}$\textit{\ }$<$\textit{\ }$0,$\textit{%
\ resp.).}

\textit{(c) Suppose }$p_{0}$\textit{\ }$=$\textit{\ }$(\alpha _{0},k_{0})$%
\textit{\ with }$k_{0}\neq 0.$\textit{\ Then there is a convex curve }$%
\gamma (s)$\textit{\ }$=$\textit{\ }$(\alpha (s),k(s)),$\textit{\ solution
to (\ref{2-1-0}) with }$\gamma (0)$\textit{\ }$=$\textit{\ }$p_{0}$\textit{\
and}%
\begin{equation*}
\lim_{s\rightarrow \pm \infty }\gamma (s)=(0,0).
\end{equation*}%
\textit{\noindent Moreover, }$\gamma $\textit{\ is symmetric with respect to
the }$k$\textit{-axis, i.e., }$(\alpha ,k)$\textit{\ }$\in $\textit{\ }$%
\gamma $\textit{\ implies }$(-\alpha ,k)$\textit{\ }$\in $\textit{\ }$\gamma
.$

\bigskip

\proof
\textbf{(of Theorem 3.3) }From (\ref{1.1}), (\ref{1.2}), we conclude that at
each point, either $k$ $=$ $0$ or%
\begin{equation}
l=-\frac{2n-i-1}{i}k.  \label{3.15}
\end{equation}

\noindent It is straightforward to see that the only stationary point for
the ODE system (\ref{2-1-0}) is $(\alpha ,k)$ $=$ $(0,0).$ A connected,
immersed umbilic hypersurface corresponding to $\alpha $ $\equiv $ $0,$ $k$ $%
\equiv $ $0$ is congruent with part of the hypersurface $E\times R$ where $E$
is a hyperplane of $R^{2n}$ by Theorem 1.5 in \cite{CCHY}.

For $\alpha (0)$ $=$ $\alpha _{0}$ $\neq $ $0,$ $k(0)$ $=$ $0,$ we have a
(unique) solution to (\ref{2-1-0}) as follows: $k(s)$ $\equiv $ $0$,%
\begin{eqnarray}
\alpha ^{\prime }(s) &=&-\alpha ^{2}(s),\text{ and hence}  \label{3.15-1} \\
\alpha (s) &=&\frac{1}{s+\alpha _{0}^{-1}}.  \notag
\end{eqnarray}

\noindent The reason is as follows. Suppose $k(s)$ $\neq $ $0$ for $\hat{s}$ 
$<$ $s$ $<$ $\hat{s}+\varepsilon $ and $k(s)$ $=$ $0$ for $0$ $\leq s\leq 
\hat{s}.$ Then (\ref{3.15}) holds for $\hat{s}$ $<$ $s$ $<$ $\hat{s}%
+\varepsilon .$ It follows by continuity that 
\begin{eqnarray*}
l(\hat{s}) &=&\lim_{s\rightarrow \hat{s}+}l(s) \\
&=&\lim_{s\rightarrow \hat{s}+}(-\frac{2n-i-1}{i}k(s)) \\
&=&-\frac{2n-i-1}{i}k(\hat{s})=0.
\end{eqnarray*}%
Now equation (\ref{2-1-0}) with $l$ given by (\ref{3.15}) holds for $\hat{s}$
$\leq $ $s$ $<$ $\hat{s}+\varepsilon .$ By the uniqueness of solutions to (%
\ref{2-1-0}) with initial data $\alpha (0)$ $=$ $\alpha _{0}$ $\neq $ $0,$ $%
k(0)$ $=$ $0,$ we get $k(s)$ $\equiv $ $0$ and $\alpha $ is given in $($\ref%
{3.15-1}$).$ We have reached a contradiction. Similarly we can argue for
negative $\hat{s}.$ If $\alpha _{0}$ $>$ $0,$ we have $\alpha (s)$ $%
\rightarrow $ $0^{+}$ as $s$ $\rightarrow $ $+\infty $ while $\alpha (s)$ $%
\rightarrow $ $+\infty $ as $s$ $\rightarrow $ ($-\alpha _{0}^{-1})^{+}.$ On
the other hand, if $\alpha _{0}$ $<$ $0,$ we have $\alpha (s)$ $\rightarrow $
$-\infty $ as $s$ $\rightarrow $ ($-\alpha _{0}^{-1})^{-}$ while $\alpha (s)$
$\rightarrow $ $0^{-}$ as $s$ $\rightarrow $ $-\infty .$ We have proved (a)
and (b).

For (c), the solution $\gamma (s)$\textit{\ }$=$\textit{\ }$(\alpha
(s),k(s)) $ to (\ref{2-1-0}) with $\gamma (0)$\textit{\ }$=$\textit{\ }$%
p_{0} $ has nonvanishing $k(s)$ for all $s$ by the uniqueness of solution to
an ODE system. By (\ref{3.15}) we can write (\ref{2-1-0}) in the following
form:%
\begin{eqnarray}
k^{\prime }(s) &=&-\frac{2n+i-1}{i}k\alpha ,  \label{3.16} \\
\alpha ^{\prime }(s) &=&\frac{2n-1}{i}k^{2}-\alpha ^{2}.  \notag
\end{eqnarray}

\noindent From (\ref{3.3}) and (\ref{3.16}) we obtain%
\begin{eqnarray}
\Pi &=&k^{2}\alpha ^{2}(\frac{2n-1}{i^{2}})(2n+3i-1)  \label{3.17} \\
&&+\alpha ^{4}(\frac{2n-i-1}{i})+(\frac{2n-1}{i}k^{2}-\alpha ^{2})^{2} 
\notag
\end{eqnarray}

\noindent since $\tilde{\sigma}_{i,n}$ $=$ $\frac{\sigma _{i,n}}{\binom{2n-2%
}{i-1}}$ $=$ $0.$ From (\ref{3.4}) and (\ref{3.15}) we have%
\begin{eqnarray}
\frac{d^{2}\alpha }{dk^{2}} &=&\frac{d}{dk}(\frac{\alpha ^{\prime }}{%
k^{\prime }})  \label{3.18} \\
&=&-\frac{\Pi }{(\frac{2n+i-1}{i}k)^{2}\alpha ^{3}}.  \notag
\end{eqnarray}

\noindent Observe that $\Pi $ $>$ $0$ by (\ref{3.17}), and hence $\frac{%
d^{2}\alpha }{dk^{2}}$ $<$ $0$ by (\ref{3.18}) for $\alpha $ $>$ $0.$ So $%
\gamma ,$ viewed as a graph $\alpha $ $=$ $\alpha (k),$ is convex (concave,
resp.) in the $\alpha $ $>$ $0$ ($\alpha $ $<$ $0,$ resp.$)$ half plane.
Therefore $\gamma $ must hit the $k$-axis. On the other hand, $\gamma $
being symmetric with respect to the $k$-axis implies that either $\gamma $
is a periodic orbit in the $k$ $>$ $0$ ($k$ $<$ $0,$ resp.$)$ plane or $%
\lim_{s\rightarrow \pm \infty }\gamma (s)=(0,0).$ In the former situation,
it ends up to have a stationary point at the positive (negative, resp.) $k$%
-axis by a topological argument. This is impossible since $(0,0)$ is the
only stationary point.

\endproof%

\bigskip

See Figure 3.3 for the $(\alpha ,k)$ diagram in the case of $c=0.$ We can
also sketch the $(\alpha ,k)$ diagram in the case of $c$ $<$ $0.$ See the
discussion in the end of Section 4.

\bigskip

\begin{figure}[h]
\includegraphics[width= 10.4cm]{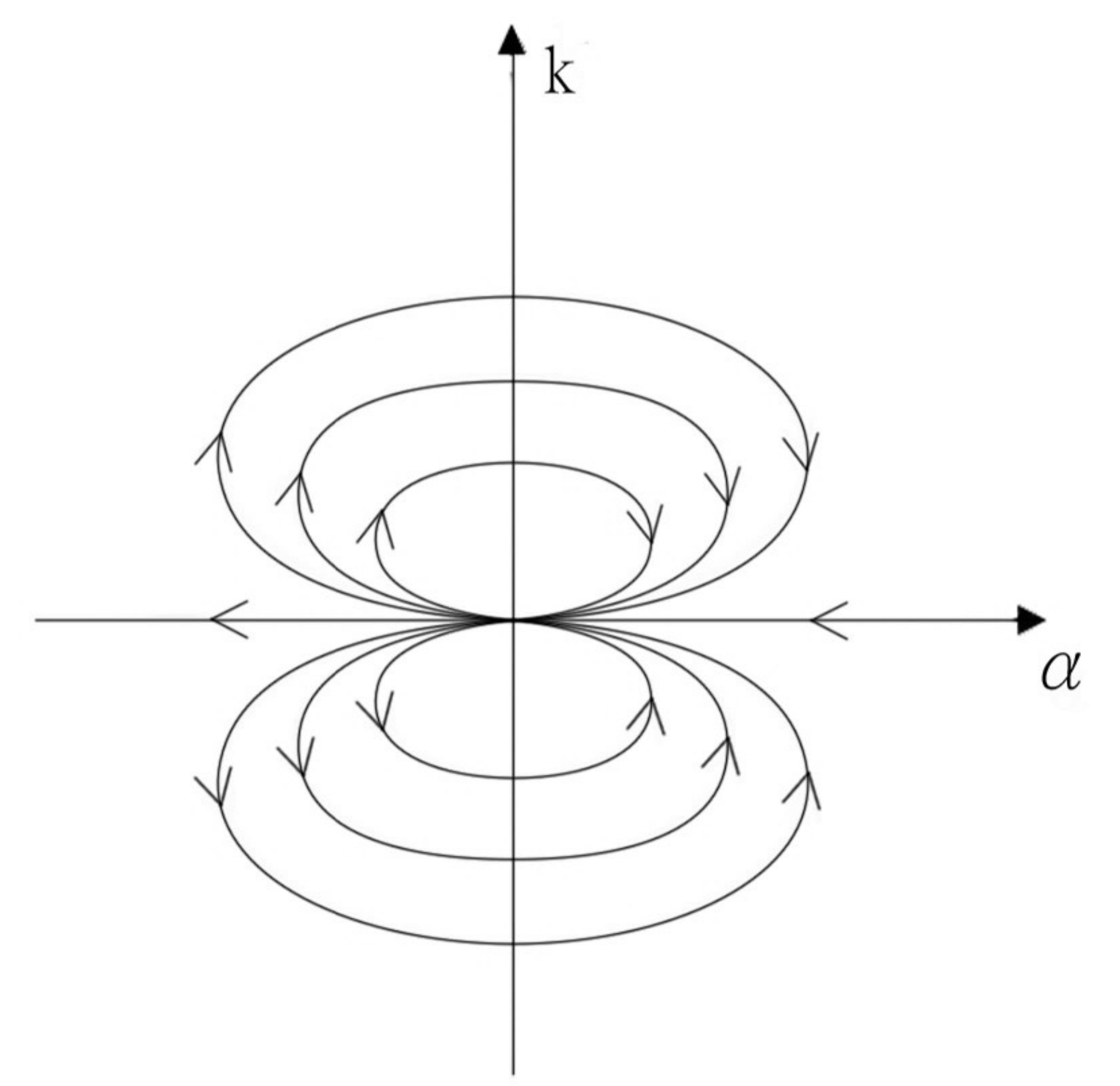} 
\caption{Figure 3.3: $c=0$, $i$ odd or even }
\end{figure}

\bigskip

\section{Proof of Theorem A and Theorem C}

The proof is in the same spirit of proving Lemma B and Theorem 1.3 (or
Theorem 1.4) in \cite{CCHY}. For completeness we will reproduce a similar
reasoning here.

\bigskip

\proof
\textbf{(of Theorem A)} For $K$ $\subset $ $\alpha k$-plane, we define the
subset $\Sigma (K)$ $\subset $ $\Sigma \backslash S_{\Sigma }$ (where $%
S_{\Sigma }$ denotes the set of singular points in $\Sigma ;$ see Section 1
for the definition of a singular point) by%
\begin{equation*}
\Sigma (K):=\{p\in \Sigma \backslash S_{\Sigma }:(\alpha (p),k(p))\in K\}.
\end{equation*}

\noindent By Proposition 4.2 in \cite{CCHY}, we learn that $k$ and $\alpha $
are constant on each leaf of the ($2n-1$)-dimensional foliation $\mathfrak{F}
$ described in Proposition 4.3 in \cite{CCHY}. On the other hand, $e_{n}$ is
transversal to the leaves by Proposition 4.3 in \cite{CCHY}, hence $\Sigma
(K)$ is open for $K$ $=$ $\{(0,0)\}$ or $\{(0,k_{c,1})\}$ for $i$ odd $%
(\{(0,0)\}$ or \{$(0,\pm k_{c,1})\},$ for $i$ even, resp.) or a periodic
orbit or an arc-type solution to (\ref{2-1-0}) in the upper half or lower
half $\alpha k$-plane or the $k$ $=$ $k_{c,2}$ line for $i$ odd ($k$ $=$ $%
\pm $ $k_{c,2}$ lines for $i$ even, resp.) by Theorem 3.1. It is also clear
that $\Sigma (K)$ is a closed set for such a $K.$

Note that $S_{\Sigma }$ consists of discrete \ (isolated singular) points by
Proposition 4.1 in \cite{CCHY}. So $\Sigma \backslash S_{\Sigma }$ is
connected and identified with $\Sigma (K)$ if $\Sigma (K)$ $\neq $ $%
\emptyset $ since $\Sigma (K)$ is open and closed.

We claim the existence of singular points. If not, $e_{n}$ is a nonvanishing
global vector field. Hence the Euler number of $\Sigma $ is zero, a
contradiction to the assumption. So $\Sigma $ contains at least a singular
point.


Observe that $\alpha \rightarrow \pm \infty $ as regular points $p_{j}$ tend
to a singular point. For $K$ being the above-mentioned sets except the $k$ $%
= $ $k_{c,2}$ line for $i$ odd ($k$ $=$ $\pm $ $k_{c,2}$ lines for $i$ even,
resp.), $\alpha $ is bounded. Therefore the only choice of $K$ for $\Sigma
(K)$ $=$ $\Sigma \backslash S_{\Sigma }$ is $K$ $=$ the $k$ $=$ $k_{c,2}$
line for $i$ odd ($k$ $=$ $\pm $ $k_{c,2}$ lines for $i$ even, resp.) (if $%
\Sigma $ has a singular point). But on such a $K,$ we have $l$ $=$ $2k.$ By
Theorem 1.3 in \cite{CCHY}, $\Sigma $ is congruent with a Pansu sphere $%
S_{\lambda }$ with $\lambda $ $=$ $k.$

\endproof%

\bigskip

We observe that if the horizontal normal $e_{2n}$ changes sign, then $\sigma
_{i,n}$ also changes sign for $i$ odd, but does not change sign for $i$ even.

\bigskip

\textbf{Lemma 4.1}. \textit{Suppose }$e_{2n}$\textit{\ changes sign (written
as }$e_{2n}\rightarrow -e_{2n}).$\textit{\ Then }$\sigma _{i,n}$\textit{\ }$%
\rightarrow -\sigma _{i,n}$\textit{\ for }$i$\textit{\ odd while }$\sigma
_{i,n}$\textit{\ }$\rightarrow \sigma _{i,n}$\textit{\ for }$i$\textit{\
even.}

\bigskip

\proof
From $e_{2n}$ $=$ $Je_{n}$ and $-\nabla _{e_{n}}e_{2n}=le_{n},$ we have $%
e_{2n}\rightarrow -e_{2n}$ implies $e_{n}\rightarrow -e_{n}$ and $%
l\rightarrow -l.$ By the definition of $\alpha :$ $\alpha e_{2n}+T\in
T\Sigma ,$ we learn that $\alpha \rightarrow -\alpha $ if $e_{2n}\rightarrow
-e_{2n}.$ It follows that the symmetric shape operator $\mathfrak{S}$ (see (%
\ref{1-1})) also changes sign. So $k\rightarrow -k.$ From formulas (\ref{1.1}%
), (\ref{1.2}), we conclude that $\sigma _{i,n}$ changes sign for $i$ odd,
but does not change sign for $i$ even.

\endproof%

\bigskip

\proof
\textbf{(of Theorem C) }To show (a), we take $K$ to be any "solution orbit"
in the ($\alpha ,k)$ diagram (see Figure 3.3). Suppose $\Sigma $ is a
closed, connected, umbilic (immersed) hypersurface of $H_{n}.$ Following the
reasoning in the proof of Theorem A, we learn that $\Sigma \backslash
S_{\Sigma }$ coincides with $\Sigma (K)$ if $\Sigma (K)$ $\neq $ $\emptyset
. $ From Theorem 3.3 (a), $K$ $=$ $\{(0,0)\}$ is impossible since $\Sigma
(K) $ in this case is unbounded while $\Sigma \backslash S_{\Sigma }$ is
bounded (noting that $\Sigma $ is compact). For $K$ being the $\alpha ^{+}$%
-axis (=$\{(\alpha ,0)$ $:$ $\alpha $ $>$ $0\}$) ($\alpha ^{-}$-axis,
resp.), $\Sigma (K)$ is also unbounded since it contains ($2n-1)$%
-dimensional spheres of radius $(\alpha ^{2}+k^{2})^{-1/2}$ = \TEXTsymbol{%
\vert}$\alpha |^{-1}$ $\rightarrow $ +$\infty $ by Theorem B. So the case
(b) in Theorem 3.3 is excluded. For case (c) in Theorem 3.3, $\Sigma (K)$ is
unbounded at two ends for the same reason since $\lim_{s\rightarrow \pm
\infty }\gamma (s)=(0,0)$ implies spheres of radius $(\alpha
^{2}+k^{2})^{-1/2}$ $\rightarrow $ $+\infty .$ We have proved (a) of Theorem
C.

To show (b), we observe that $\sigma _{i,n}$ changes sign for $i$ odd if we
change the sign of $e_{2n}$ while $l$ and $k$ change sign, and $\alpha $
does not change. Moreover, when $i$ is odd, $(\alpha (s),k(s))$ is a
solution to (\ref{2-1-0}) for $\sigma _{i,n}$ $=$ $c$ $<$ $0$ if and only if 
$(\alpha (s),-k(s))$ is a solution to (\ref{2-1-0}) for $\sigma _{i,n}$ $=$ $%
c$ $>$ $0.$ So a similar reasoning as in the proof of Theorem A shows that
only $K$ = $\{k=-k_{-c,1}\}$ is a possible "solution orbit" and $\Sigma $%
\textit{\ }must be a Pansu sphere up to a Heisenberg translation.

\endproof%

\bigskip

For $\sigma _{i,n}$ $=$ $c$ $<$ $0$ and $i$ even, we discuss and sketch the $%
(\alpha ,k)$-diagram as follows. First observe that when $c$ $<$ $0,$ $k$ is
never zero for $i$ $\geq $ $2.$ From (\ref{1.1}) we get%
\begin{equation}
l=\frac{c}{\binom{2n-2}{i-1}}k^{1-i}-\frac{2n-i-1}{i}k  \label{4.1}
\end{equation}

\noindent for $2\leq i\leq 2n-1$. By (\ref{4.1}) we compute%
\begin{eqnarray}
l-2k &=&k[\frac{c}{\binom{2n-2}{i-1}}k^{-i}-\frac{2n-i-1}{i}-2],  \label{4.2}
\\
k(k-l) &=&k^{2}[1-\frac{c}{\binom{2n-2}{i-1}}k^{-i}+\frac{2n-i-1}{i}]. 
\notag
\end{eqnarray}

\noindent Since $c$ $<$ $0$ and (hence) $k$ $\neq $ $0$, we have $k(k-l)$ $>$
$0$ by the second equality of (\ref{4.2}). It follows that there are no
stationary points (where $k^{\prime }$ $=$ $0,$ $\alpha ^{\prime }$ $=$ $0)$
for the ODE system (\ref{2-1-0}) in this situation ($c$ $<$ $0$ and $i$
even). Recall that $\Psi _{0}^{\pm }$ in $(\alpha ,k)$-plane is defined by $%
k^{2}-\alpha ^{2}-kl$ $=$ $0$ (where $\alpha ^{\prime }$ $=$ $0$), in $%
k\gtrless $ 0 regions, resp.. For the case $i$ $=$ $2,$ $\Psi _{0}^{\pm }$
is described by the following hyperbolic curve:%
\begin{equation*}
(1+\frac{2n-i-1}{i})k^{2}-\alpha ^{2}=\frac{c}{\binom{2n-2}{i-1}}.
\end{equation*}

\noindent On the other hand, for $i$ $\geq $ $4,$ $\Psi _{0}^{\pm }$ is a
curve of different type. We have sketched the solution diagram in Figure 4.1
and Figure 4.2 for these two cases.

\bigskip

\begin{figure}[h]
\includegraphics[width= 10.4cm]{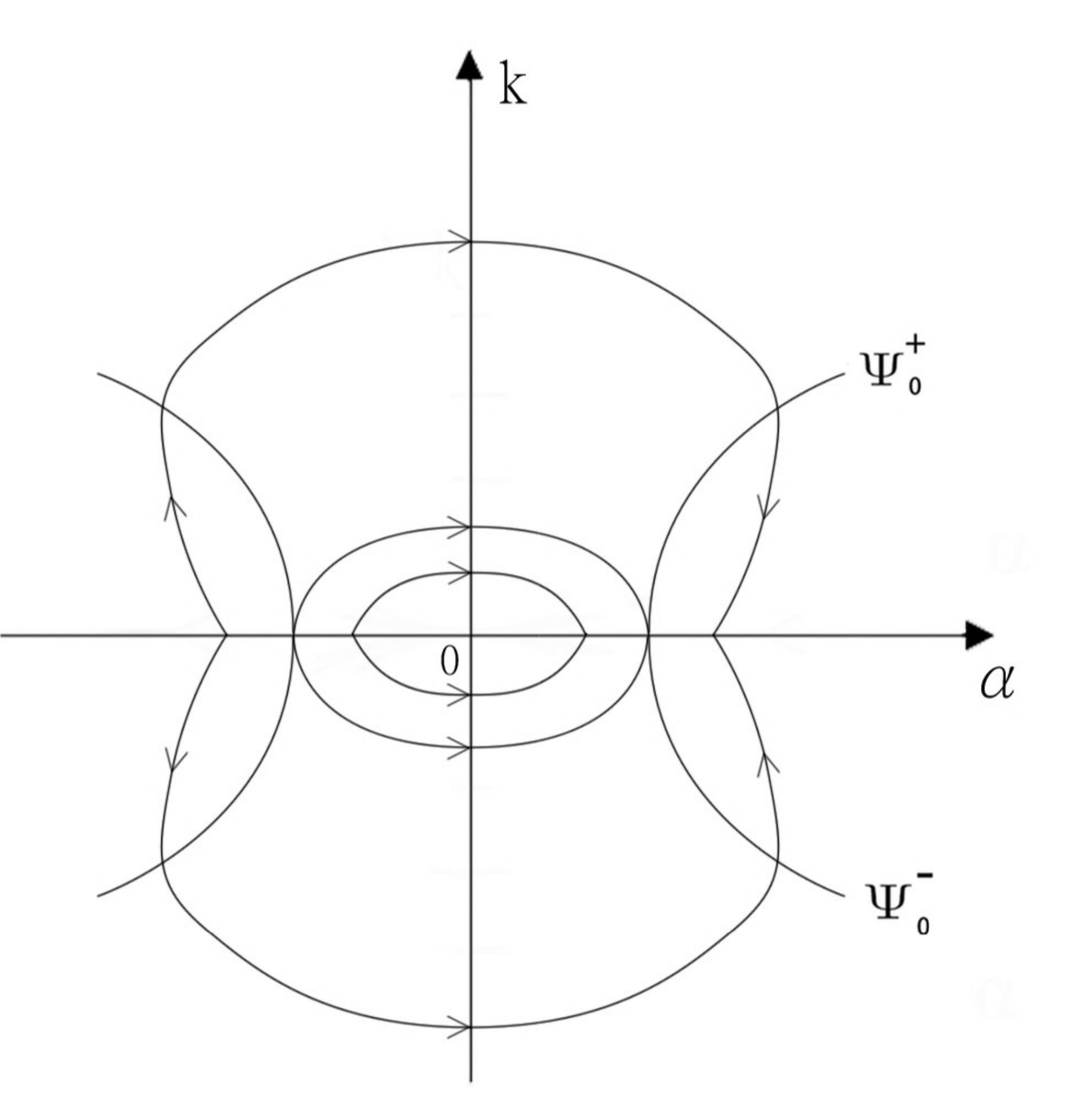} 
\caption{Figure 4.1: $c<0$, $i=2$}
\end{figure}

\begin{figure}[h]
\includegraphics[width= 10.4cm]{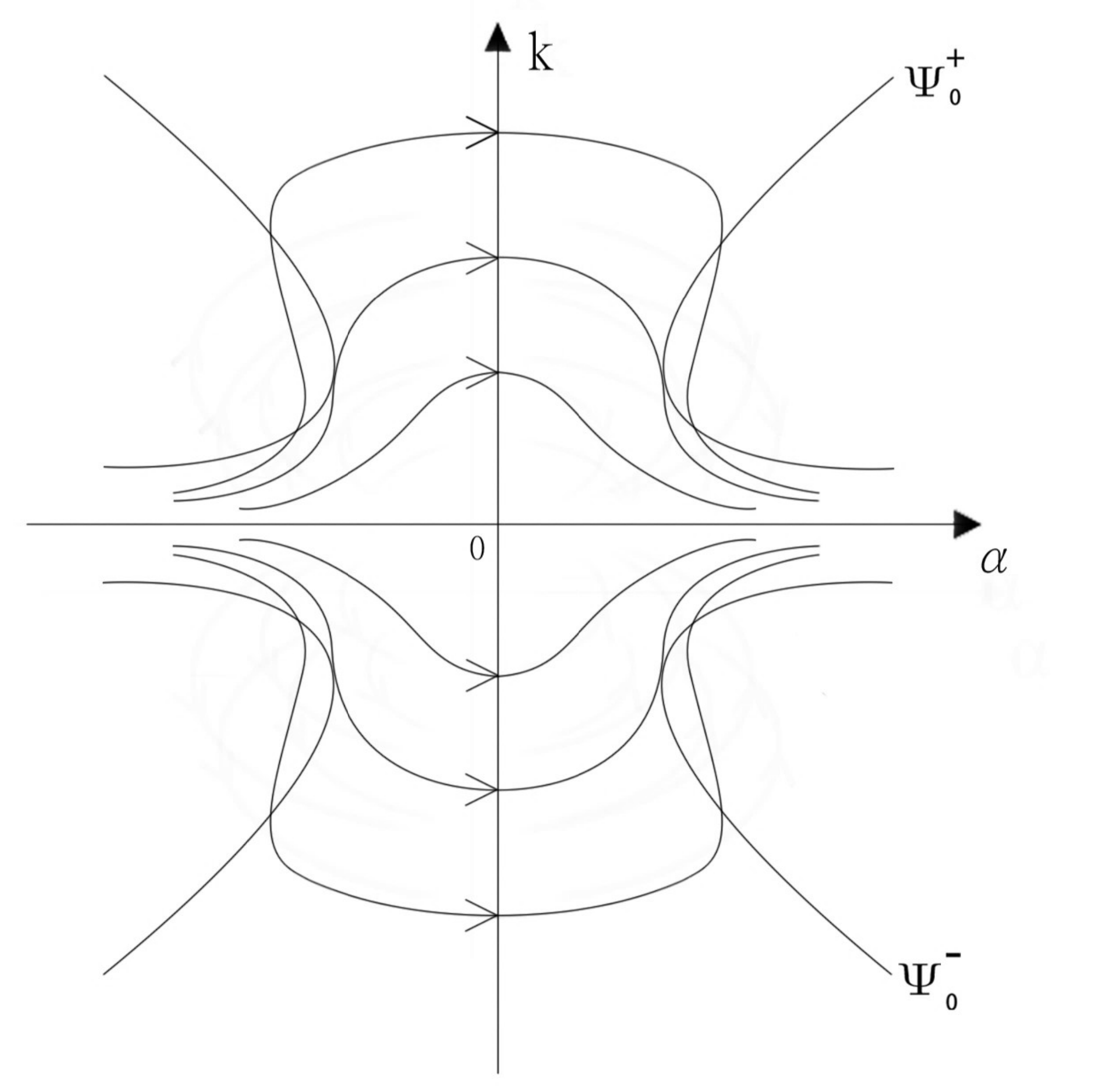} 
\caption{Figure 4.2: $c<0$, even $i \ge 4$}
\end{figure}

\bigskip

We say a hypersurface $\Sigma $ can be completed to a closed hypersurface $%
\bar{\Sigma}$ if $\bar{\Sigma}$ contains $\Sigma $ and $\bar{\Sigma}$ is
closed, i.e., compact with no boundary.

\bigskip

\textbf{Proposition 4.2.} \textit{Given a solution curve }$\gamma $\textit{\
in Figure 4.2: }$c<0,$\textit{\ even }$i\geq 4$ \textit{(note that }$k$%
\textit{\ }$\equiv $\textit{\ }$0$\textit{\ (}$\alpha $\textit{-axis) is not
a solution curve)}$,$\textit{\ suppose there corresponds an umbilic
hypersurface }$\Sigma $ \textit{with its }$(\alpha (s),k(s))$\textit{\ curve
being }$\gamma .$ \textit{Then }$\Sigma $\textit{\ can be completed to a} 
\textit{closed, rotationally invariant hypersurface }$\bar{\Sigma}$ \textit{%
up to a Heisenberg translation and }$\bar{\Sigma}$\textit{\ is }$C^{2}$%
\textit{\ smooth at its two singular points.}

\bigskip

\proof
By Theorem B (c), $\Sigma $ is congruent with part of a rotationally
invariant hypersurface $\Sigma ^{\prime }$ and the radius $r$ $=$ $\frac{1}{%
\sqrt{\alpha ^{2}+k^{2}}}$ of leaves in the associated foliation strictly
decreases to $0$ as $\alpha $ goes to $\pm \infty $ (while $s$ $\rightarrow $
$\pm \infty )$. Write $\Sigma ^{\prime }$ as a graph $t$ $=$ $g_{+}(r)$ ($%
g_{-}(r),$ resp.)$,$ $g_{\pm }$ $\in $ $C^{\infty },$ for $r$ $:=$ $[\Sigma
_{j=1}^{n}(x_{j}^{2}+y_{j}^{2})]^{1/2}$ close to $0$ and $\alpha $ near $%
+\infty $ ($\alpha $ near $-\infty ,$ resp.)$.$ For $\vec{F}$ $=$ $-y_{j}%
\frac{\partial }{\partial x_{j}}$ $+$ $x_{j}\frac{\partial }{\partial y_{j}}%
, $ we compute 
\begin{equation*}
\nabla g_{\pm }+\vec{F}=(g_{\pm }^{\prime }(r)\frac{x_{j}}{r}-y_{j})\frac{%
\partial }{\partial x_{j}}+(g_{\pm }^{\prime }(r)\frac{y_{j}}{r}+x_{j})\frac{%
\partial }{\partial y_{j}}.
\end{equation*}

\noindent It follows that 
\begin{equation}
\frac{1}{|\alpha |}=|\nabla g_{\pm }+\vec{F}|=\sqrt{(g_{\pm }^{\prime
}(r))^{2}+r^{2}}\rightarrow 0  \label{4-3}
\end{equation}

\noindent as $s$ $\rightarrow $ $\pm \infty .$ So $r$ $\rightarrow $ $0$ and 
$g_{\pm }^{\prime }(r)$ $\rightarrow $ $0.$ Therefore $t$ is bounded as $r$ $%
\rightarrow $ $0$ and $\lim_{r\rightarrow 0}g_{\pm }(r)$ exists, denoted as $%
g_{\pm }(0)$. Also $g_{\pm }$ is $C^{1}$ smooth at $r$ $=$ $0$ (singular
points)$.$ So $\Sigma ^{\prime }$ can be completed to a closed rotationally
invariant hypersurface $\bar{\Sigma}$ containing two (singular) points $(r$ $%
=$ $0,g_{\pm }(0))$.

From the formula in (\ref{4-3}) and $r$ $=$ $\frac{1}{\sqrt{\alpha ^{2}+k^{2}%
}}$, we compute%
\begin{eqnarray}
\left( \frac{g_{\pm }^{\prime }(r)}{r}\right) ^{2}+1 &=&\frac{1}{\alpha
^{2}r^{2}}  \label{4-4} \\
&=&\frac{\alpha ^{2}+k^{2}}{\alpha ^{2}}\rightarrow 1  \notag
\end{eqnarray}

\noindent as $\alpha \rightarrow \pm \infty ,$ $k\rightarrow 0$ (while $s$ $%
\rightarrow $ $\pm \infty ).$ It follows from (\ref{4-4}) that $\frac{g_{\pm
}^{\prime }(r)}{r}\rightarrow 0$ as $r\rightarrow 0$. Hence $g_{\pm
}^{^{\prime \prime }}(0)$ $=$ $0.$ A routine verification shows that $\bar{%
\Sigma}$ defined by $g_{\pm }$ is $C^{2}$ smooth at two (singular) points $%
(r $ $=$ $0,g_{\pm }(0))$, resp..

\endproof%

\bigskip

We remark that the above hypersurface $\bar{\Sigma}$ in Proposition 4.2 is
not a Pansu sphere on which $k$ is a constant.

For the case $c$ $<$ $0$ and odd $i$ $\geq $ $3,$ the $(\alpha ,k)$ diagram
is the reflection of Figure 3.1 with respect to the $\alpha $-axis by Lemma
4.1 (see Figure 4.3).

\bigskip

\begin{figure}[h]
\includegraphics[width= 10.4cm]{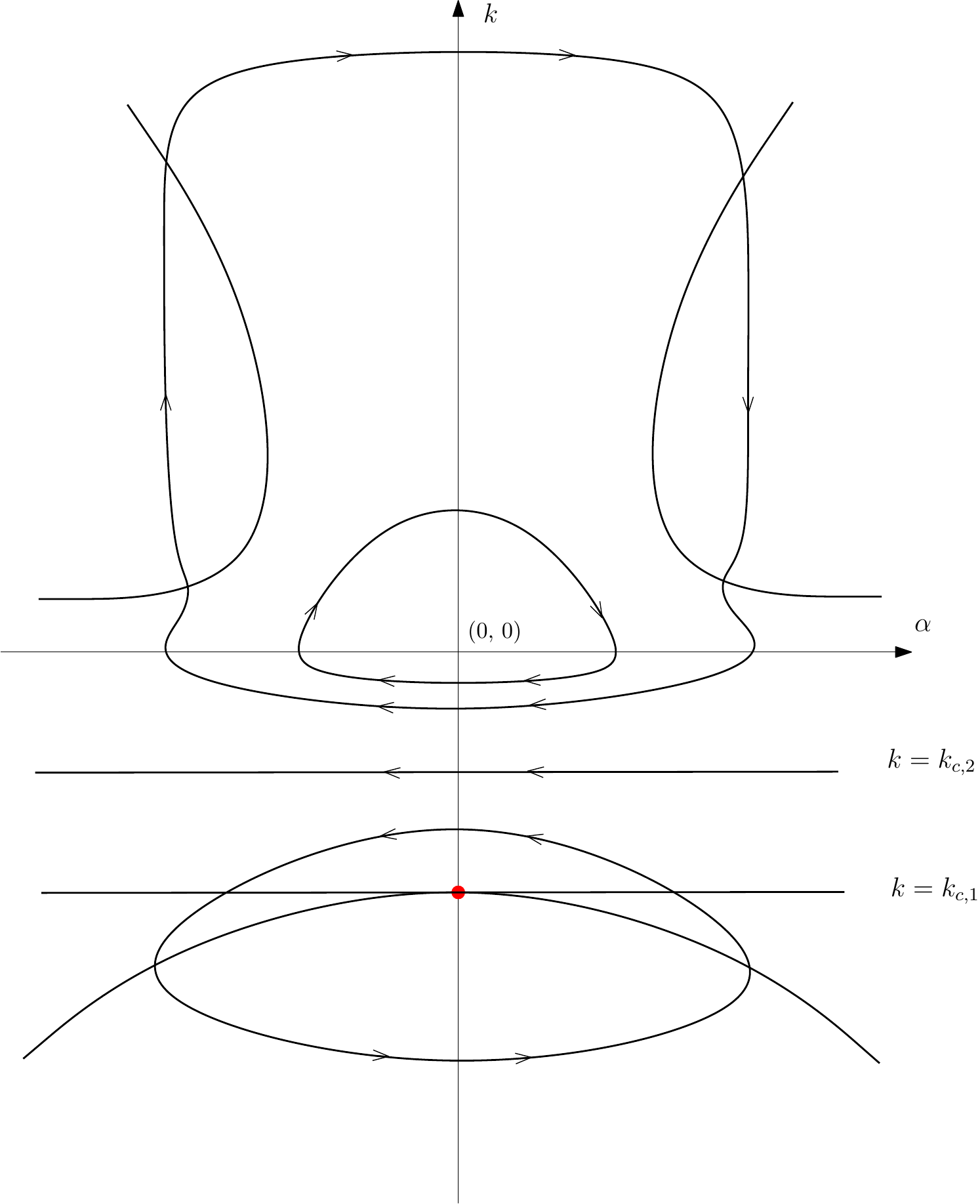} 
\caption{Figure 4.3: $c<0$, odd $i\geq 3$}
\end{figure}

\bigskip

\section{Appendix A: closed form of solutions to (\protect\ref{2-1-0})}

Write (\ref{2-1-0}) in terms of differential form as follows:%
\begin{equation}
(k^{2}-kl-\alpha ^{2})dk+(2k-l)\alpha d\alpha =0  \label{B1}
\end{equation}

\noindent where we recall (see (\ref{2-1-1})) that 
\begin{equation}
l=\frac{c}{\binom{2n-2}{i-1}k^{i-1}}-\frac{2n-i-1}{i}k.  \label{B1-0}
\end{equation}

\noindent Try to find an integrating factor $f(k)$ for (\ref{B1}). Let $%
\zeta $ $:=$ $\alpha ^{2}.$ That is, there holds%
\begin{equation}
d\{f(k)[(k^{2}-kl-\zeta )dk+(2k-l)\frac{d\zeta }{2}]\}=0.  \label{B1-1}
\end{equation}

\noindent It follows that 
\begin{equation*}
\frac{\partial (f(k)\zeta )}{\partial \zeta }+\frac{\partial }{\partial k}%
(f(k)\frac{2k-l}{2})=0.
\end{equation*}

\noindent So $f(k)$ satisfies the following ODE:%
\begin{equation}
2f(k)+\frac{d}{dk}(f(k)(2k-l))=0.  \label{B2}
\end{equation}

\noindent Let $g(k)$ $:=$ $f(k)(2k-l).$ We can rewrite (\ref{B2}) as%
\begin{equation*}
\frac{dg(k)}{dk}=-\frac{2g(k)}{2k-l},
\end{equation*}

\noindent and get the solution:%
\begin{equation}
g(k)=\exp \{-\int \frac{2}{2k-l}dk\}.  \label{B3}
\end{equation}

\noindent From (\ref{B1-1}) and (\ref{B2}), we obtain%
\begin{equation*}
d(\int f(k)(k^{2}-kl)dk+\frac{\zeta g(k)}{2})=0
\end{equation*}

\noindent Thus we have%
\begin{equation*}
\int f(k)(k^{2}-kl)dk+\frac{\alpha ^{2}g(k)}{2}=C_{1}
\end{equation*}

\noindent for some constant $C_{1}.$ We can then express $\alpha $ in terms
of $k$ as follows:%
\begin{equation}
\alpha ^{2}=\frac{2C_{1}}{g(k)}-\frac{2}{g(k)}\int f(k)(k^{2}-kl)dk
\label{B4}
\end{equation}
\noindent Substituting (\ref{B1-0}) into (\ref{B3}) gives%
\begin{eqnarray*}
g(k) &=&\exp \{\frac{-2}{i}\int \frac{d(k^{i})}{(\frac{2n-i-1}{i})k^{i}-%
\frac{c}{\binom{2n-2}{i-1}}}\} \\
&=&C_{2}\exp \{\frac{-2}{2n+i-1}\log \text{ }|k^{i}-\frac{i}{2n+i-1}\frac{c}{%
\binom{2n-2}{i-1}}|
\end{eqnarray*}

\noindent for some constant $C_{2}$ $>$ $0,$ and hence we have%
\begin{equation}
g(k)=C_{2}|k^{i}-\frac{i}{2n+i-1}\frac{c}{\binom{2n-2}{i-1}}|^{\frac{-2}{%
2n+i-1}}\text{.}  \label{B5}
\end{equation}

\noindent Next by (\ref{B5}) we compute%
\begin{eqnarray}
&&\int f(k)(k^{2}-kl)dk  \label{B6} \\
&=&\int g(k)\frac{k^{2}-kl}{2k-l}dk  \notag \\
&=&C_{2}\int \{|k^{i}-\frac{i}{2n+i-1}\frac{c}{\binom{2n-2}{i-1}}|^{\frac{-2%
}{2n+i-1}}(k^{i}-\frac{i}{2n+i-1}\frac{c}{\binom{2n-2}{i-1}})^{-1}  \notag \\
&&\frac{k}{2n+i-1}[(2n-1)k^{i}-\frac{c}{\binom{2n-2}{i-1}}i]\}dk.  \notag
\end{eqnarray}

\noindent In view of (\ref{B5}) and (\ref{B6}), (\ref{B4}) is a closed form
of solutions $(\alpha (s),k(s))$ to (\ref{2-1-0}).

\end{document}